\documentclass{article}
\usepackage{amsmath}
\usepackage{eurosym}
\usepackage{makeidx}
\usepackage{amssymb}

\newtheorem{theorem}{Theorem}

\newtheorem{corollary}[theorem]{Corollary}

\newtheorem{definition}[theorem]{Definition}

\newtheorem{lemma}[theorem]{Lemma}

\newtheorem{proposition}[theorem]{Proposition}
\newtheorem{remark}[theorem]{Remark}

\input{tcilatex}
\begin{document}

\title{Numbers and numerosities}
\author{Vieri Benci\thanks{%
Dipartimento di Matematica, Universit\`{a} degli Studi di Pisa, Via F.
Buonarroti 1/c, 56127 Pisa, Italy, e-mail: vieri.benci@unipi.it.}}
\maketitle

\begin{abstract}
We develop new aspects of the the of numerosity theory; more exactly, we
emphasize its relation with the ordinal numbers, cardinal numbers, hyperreal
numbers and surreal numbers. In particular, we combine the notion of
numerosity with the idea of continuum and we get a definition of Euclidean
line which includes all the sets of infinite numbers mentioned above.

\textbf{keywords: }Euclid's principle, Hume's principle, numerosities,
ordinal numbers, cardinal numbers, hyperreal numbers, surreal numbers,
Euclidean line.
\end{abstract}

\tableofcontents

\section{Introduction}

The theory of numerosity, introduced in \cite{benci95b,BDN2003} is a
development of the classical cardinality for measuring the size of infinite
sets. Unlike Cantor's theory, where all countably infinite sets have the
same cardinality, numerosity theory allows a finer distinction between
different infinite sets and it preserves properties more aligned with our
intuitive understanding of "number of elements." The notion of numerosity
has been developed in several direction; see e.g. \cite%
{BDN2018,BDNF1,BF,BLDN,DNFtup,DNF,BLL,mancosu09,mancosu16,mancosu10,mancosu23}
and the references therein.

In this paper we summarize some old results and present new aspects. In the
first part we introduce the theory of numerosity in a new way. First we
analyze the idea of number and we define the structure of \textit{counting
system}; in this context, we present the notion of \textit{numerosity} as a
natural extension of the notion of natural number. In particular, we
emphasize its relation with the ordinal numbers, cardinal numbers, hyperreal
numbers and surreal numbers. The relation between numerosities, hyperreal
numbers and surreal numbers is a new fulfillment and suggests further
developments. Then, we build a model of the numerosities based on labelled
sets as in \cite{benci95b,BDN2003,BLL} which allows to develop new features.

\subsection{Preliminaries ad notation\label{pn}}

In principle it would be desirable to define the operation of "counting" for
the class of all sets; however, in order to develop the theory, it is
convenient to work in a "universe" which is itself a set provided that it is
very large. Hence we will work in a universe in which the classes are
replaced by sets of strongly inaccessible cardinality.

\begin{definition}
\label{uno}A cardinal number $\chi $ is \textbf{inaccessible} if it is not a
sum of fewer than $\chi $ cardinals that are less than $\chi $ and $\zeta
<\chi $ implies $2^{\zeta }<\chi $.\ $\chi $ is \textbf{strongly
inaccessible\ }if it is inaccessible and uncountable.
\end{definition}

The first inaccessible cardinal number is $\aleph _{0}.$ The first strongly
inaccessible cardinal number will be denoted by $\mathbf{\kappa }.$ The
existence of sets of strongly inaccessible cardinality is established by the
Axiom of Inaccessibility which is independent from ZFC. We will assume this
axiom and, in this paper, a set having strongly inaccessible cardinality
will be called \textbf{pseudoclass}.

We will assume that there exists a set of atoms\footnote{%
In set theory, an atom $a$ is any entity that is not a set, namely $a$ is an
atom if and only if%
\begin{equation*}
\forall x,\ x\notin a
\end{equation*}%
} $\mathbf{Ato}$ having cardinality $\mathbf{\kappa }$. Moreover we will
denote by $\mathbf{Card\ }$the pseudoclass of the first $\mathbf{\kappa }$
cardinal numbers an by $\mathbf{Ord}$ the pseudoclass of the first \textbf{$%
\kappa $} ordinal numbers. The cardinality of a set $X$ will be denoted by $%
\left\vert X\right\vert $.

Given any set $E,$ let $V_{\beta }(E),\ \beta \in \mathbf{Ord,}$ be the
superstructure on $E$, namely the family of sets which is inductively
defined as follows:%
\begin{eqnarray*}
V_{0}(E) &=&E; \\
V_{\beta +1}(E) &=&V_{\beta }(E)\cup \mathbf{\wp }\left( V_{\beta
}(E)\right) ; \\
V_{\beta }(E) &=&\bigcup\limits_{\gamma <\beta }V_{\gamma }(E)\ \text{if }%
\beta \text{ is a limit ordinal.}
\end{eqnarray*}%
If $A\in V_{\beta }(\mathbf{Ato})\backslash \bigcup\limits_{\gamma <\beta
}V_{\gamma }(\mathbf{Ato})$, we say that $A$ has rank $\beta $ and we will
write $rank(A)=\beta .$ Now we set:%
\begin{equation*}
\mathbf{U}_{\max }:=\bigcup\limits_{\gamma <\beta }\left\{ E\in V_{\gamma }(%
\mathbf{Ato})\ |\ \ \left\vert E\right\vert <\mathbf{\kappa }\right\}
\end{equation*}%
It is well known that $\mathbf{U}_{\max }$ is a model of $ZFC$ and we can
work in it without the problems related to the theory of classes.

\begin{definition}
A set $\mathbf{U}\subseteq \mathbf{U}_{\max }\backslash \mathbf{Ato}$ is
called universe (of sets) if

\begin{enumerate}
\item $A,B\in \mathbf{U}\Rightarrow A\cup B\in \mathbf{U}$

\item $A,B\in \mathbf{U}\Rightarrow A\times B\in \mathbf{U}$

\item $A\in \mathbf{U\ }$and$\mathbf{\ \ }B\subseteq A\Rightarrow B\in 
\mathbf{U}$
\end{enumerate}
\end{definition}

Moreover we shall use the following notation:

\begin{itemize}
\item $\mathbf{Inf}$ will denote the family of infinite subsets in $\mathbf{U%
}_{\max }.$

\item $\mathbf{Fin}$\textbf{\ }will denote the family of the finite subsets
in $\mathbf{U}_{\max }.$

\item if $E\in \mathbf{U,\ }\mathfrak{\wp }(E)$ will denote the family of
subsets of $E.$

\item $\mathfrak{\wp }_{\omega }(E):=\mathfrak{\wp }(E)\cap \mathbf{Fin}$
will denote the family of finite subsets of $E$

\item $\mathfrak{F}(A,B)$ will denote the family of functions $%
f:A\rightarrow B,$ namely $\mathfrak{F}(A,B)=B^{A}.$
\end{itemize}

\section{Counting systems\label{CS}}

\subsection{The numbers}

One of the aims in counting the elements of sets is the comparison of their
sizes. We denote by $\preceq $ is a total preorder\footnote{%
We recall that a total preorder is a\textit{\ trnsitive} and\textit{\
reflexive} relation such that $x\preceq y\ $or $y\preceq x$.} relation over $%
\mathbf{U}$ and, as usual, we set%
\begin{equation*}
X\cong Y:\Leftrightarrow \left( X\preceq Y\ \ \text{and\ \ }Y\preceq
X\right) .
\end{equation*}

\begin{definition}
\label{C}A \textbf{comparison system} is a couple $(\mathbf{U},\preceq )$
where $\mathbf{U}$ is a universe and $\preceq $ is a total preorder
relation, called \textbf{comparison relation}, which satisfies the following
properties:

\begin{enumerate}
\item \label{C1}\textsf{Null Set Principle:} $A\cong \varnothing ,$ if and
only if $A=\varnothing .$

\item \label{C3}\textsf{Union Principle:}\ If $A\cap B=A^{\prime }\cap
B^{\prime }=\varnothing ,\ $and $A\cong A^{\prime },\ B\cong B^{\prime },$
then%
\begin{equation*}
A\cup B\cong A^{\prime }\cup B^{\prime }
\end{equation*}

\item \label{C4}\textsf{Product Principle:}\ If $A\cong A^{\prime }$ and $%
B\cong B^{\prime }$;then%
\begin{equation*}
A\times B\cong B\times A\cong A^{\prime }\times B^{\prime }.
\end{equation*}

\item \label{C5}\textsf{Unit Principle: }%
\begin{equation*}
\{a\}\times B\cong B
\end{equation*}

\item \label{C2}\textsf{Comparison Principle: }$A\preceq B,$ if and only if
there exists a set $A^{\prime }\subseteq B$ such that 
\begin{equation*}
A\cong A^{\prime }.
\end{equation*}
\end{enumerate}
\end{definition}

If we have a comparison system, then it is possible to build the notion of
number:

\begin{definition}
\label{SN}A set of numbers $\mathcal{N}$ is a set of atoms such that there
exists a bijection 
\begin{equation*}
\Phi :\mathbf{U}/\cong \ \rightarrow \mathcal{N}
\end{equation*}
\end{definition}

Then given a set $A,$ the number of its elements is given by 
\begin{equation*}
\mathfrak{n}(A)=\Phi \left( \left[ A\right] _{\cong }\right)
\end{equation*}%
Notice that in every set of numbers there are two distinguished elements:%
\begin{equation*}
0:=\mathfrak{n}(\varnothing )
\end{equation*}%
and 
\begin{equation*}
1:=\mathfrak{n}(\left\{ \varnothing \right\} ).
\end{equation*}%
Since we have required $\preceq $ to be a total preorder relation, the
following fact follows:

\begin{proposition}
The numbers is a linearly ordered set with respect to the following order
relation: given $\alpha =\mathfrak{n}(A)$ and $\beta =\mathfrak{n}(B)$ 
\begin{equation*}
\alpha \leq \beta :\Leftrightarrow A\preceq B.
\end{equation*}
\end{proposition}

\subsubsection{Operations with numbers\label{ON}}

Given a set of numbers we can define also the two basic operations: the sum
and the product.

\begin{definition}
\label{S}Given two numbers $\alpha =\mathfrak{n}(A)$ and $\beta =\mathfrak{n}%
(B)\ $with\ $A\cap B=\varnothing $, we set 
\begin{equation*}
\alpha +\beta :=\mathfrak{n}(A\cup B)
\end{equation*}
\end{definition}

By the Union Principle, it follows that the operation "+" is well posed,
namely it is independent of the choice of $A$ and $B$. It is immediate to
check that the sum satisfies the commutative property and the associative
property since%
\begin{equation*}
A\cup B=B\cup A
\end{equation*}%
and%
\begin{equation*}
\left( A\cup B\right) \cup C=A\cup \left( B\cup C\right)
\end{equation*}%
Moreover we have that $0=\mathfrak{n}(\varnothing )$ is the identity element
for the sum. It is obvious that two numbers can always be added and we have
that $\alpha +\beta =\mathfrak{n}(A)+\mathfrak{n}(B)$: in fact if $A\cap
B\neq \varnothing ,$ we can replace $A$ with the set $A\times \{c\}$ with $%
c\notin A\cup B$;\ then, $\left( A\times \{c\}\right) \cap B=\varnothing ,$
and hence,%
\begin{equation*}
\alpha +\beta =\mathfrak{n}(\left( A\times \{c\}\right) \cup B)=\mathfrak{n}%
(A\times \{c\})+\mathfrak{n}(B)=\mathfrak{n}(A)+\mathfrak{n}(B)
\end{equation*}

\begin{definition}
\label{P}Given two numbers $\alpha =\mathfrak{n}(A)$ and $\beta =\mathfrak{n}%
(B)\ $with\ $A,B$ as in Def. \ref{C}-(\ref{C4}), we set 
\begin{equation*}
\alpha \cdot \beta :=\mathfrak{n}(A\times B)
\end{equation*}
\end{definition}

By the Union Principle, it follows that also the operation $"\cdot "$ is
well posed. In Def. \ref{C}, we have imposed the commutative property. The
associative property is satisfied if we identify the couple $((a,b),c)$ with 
$(a,(b,c)).$ This property will be fulfilled by all our examples.

In conclusion, the numbers as defined by Def. \ref{SN}, satisfy the basic
algebraic properties (BAC):

\begin{itemize}
\item commutative property with respect to + and $\cdot $

\item associative property with respect to + and $\cdot $

\item existence of the identity elements "0" and "1" with respect to + and $%
\cdot $

\item distributive property
\end{itemize}

\begin{definition}
\label{CT}The triple $(\mathbf{U},\mathcal{N},\mathfrak{n})$ constructed
over a comparison relation, is called \textbf{counting system.}
\end{definition}

\subsection{The main counting systems}

\subsubsection{The finite counting systems}

Now let us see some examples of counting systems:

\bigskip

\textbf{Example 1 -} we take

\begin{itemize}
\item $\mathbf{U}=\mathbf{U}_{\max }$:

\item $\mathcal{N}=\{0,1,2,M\}$ where the number $M$ is read \textbf{"many".}
\end{itemize}

Then, there exists a unique comparison relation which induces the following
arithmetic:

\begin{equation*}
\begin{tabular}{|l|l|l|l|l|}
\hline
$+$ & $\left[ 0\right] $ & $\left[ 1\right] $ & $\left[ 2\right] $ & $\left[
M\right] $ \\ \hline
$\left[ 0\right] $ & $0$ & $1$ & $2$ & $M$ \\ \hline
$\left[ 1\right] $ & $1$ & $2$ & $M$ & $M$ \\ \hline
$\left[ 2\right] $ & $2$ & $M$ & $M$ & $M$ \\ \hline
$\left[ M\right] $ & $M$ & $M$ & $M$ & $M$ \\ \hline
\end{tabular}%
;\ \ 
\begin{tabular}{|l|l|l|l|l|}
\hline
$\cdot $ & $\left[ 0\right] $ & $\left[ 1\right] $ & $\left[ 2\right] $ & $%
\left[ M\right] $ \\ \hline
$\left[ 0\right] $ & $0$ & $0$ & $0$ & $0$ \\ \hline
$\left[ 1\right] $ & $0$ & $1$ & $2$ & $M$ \\ \hline
$\left[ 2\right] $ & $0$ & $2$ & $M$ & $M$ \\ \hline
$\left[ M\right] $ & $0$ & $M$ & $M$ & $M$ \\ \hline
\end{tabular}%
\end{equation*}

Actually, this example, does not provide the "smallest" counting system
since we can take $\mathcal{N}=\left\{ 0,1\right\} .$ In this case, $(%
\mathbf{U},\left\{ 0,1\right\} ,\mathfrak{n})$ reduces to a Boolean algebra
where, in this case, the addition "+" corresponds to \textit{"or"} ($\vee $)
and the product "$\cdot $" corresponds to to \textit{"and"} ($\wedge $):%
\begin{equation}
\begin{tabular}{|l|l|l|}
\hline
$+\equiv \vee $ & $\left[ 0\right] $ & $\left[ 1\right] $ \\ \hline
$\left[ 0\right] $ & $0$ & $1$ \\ \hline
$\left[ 1\right] $ & $1$ & $1$ \\ \hline
\end{tabular}%
;\ \ 
\begin{tabular}{|l|l|l|}
\hline
$\cdot \equiv \wedge $ & $\left[ 0\right] $ & $\left[ 1\right] $ \\ \hline
$\left[ 0\right] $ & $0$ & $0$ \\ \hline
$\left[ 1\right] $ & $0$ & $1$ \\ \hline
\end{tabular}%
.  \label{du}
\end{equation}

\subsubsection{The Euclid's Principle and the Natural numbers}

If we want to exclude these interesting, but mathematically trivial
examples, we need to add some other principle: for example the \textit{V
common notion} of Euclid's elements:

\begin{center}
\textit{The whole is greater than the part.}
\end{center}

In our formalism

\begin{itemize}
\item \textbf{Euclid's principle} - \textit{Given two sets }$F$\textit{\ and 
}$G$\textit{\ such that }$F$\textit{\ is a proper part of }$G$\textit{, then}
$F\prec G$.
\end{itemize}

The most important counting system which satisfies Euclid's Principle is the
counting system of natural numbers $(\mathbf{Fin},\mathbb{N},\left\vert
\cdot \right\vert )$ where

\begin{itemize}
\item $\mathbf{Fin}$ is the family of finite sets.

\item $\mathbb{N}$ is the set of natural numbers.

\item $\left\vert A\right\vert =\mathfrak{n}\left( A\right) $ is the number
of elements of a finite set.
\end{itemize}

$(\mathbf{Fin},\mathbb{N},\left\vert \cdot \right\vert )$ is the smallest
counting system which satisfies the Euclid's Principle namely every counting
system which satisfies the Euclid's Principle contains $(\mathbf{Fin},%
\mathbb{N},\left\vert \cdot \right\vert )$ as a subsystem.

The system of natural numbers satisfies a good algebra; hence, we can easily
buid $\mathbb{Z}$, $\mathbb{Q}$ and $\mathbb{R}$ which satsfy the very rich
algebra which is used by everybody.

\subsubsection{The Euler's infinite}

In order to get a theory which incudes infinite sets, Euler introduced the
symbol "$\infty $" which is similar to the number "$M$" in example 1. Using
Euler's infinite, we obtain the counting system $(\mathbf{U}_{\max },\mathbb{%
N}\cup \left\{ \infty \right\} ,\mathfrak{n})$ with the relations%
\begin{eqnarray*}
n+\infty &=&\infty ,\  \\
0\cdot \infty &=&0, \\
n\cdot \infty &=&\infty \ \text{for}\ n\neq 0
\end{eqnarray*}%
This system, satisfies the Euclid's pronciple only for finite sets, and
hence it does not have a good algebra: in particular, the equation%
\begin{equation*}
x+\infty =\infty ,
\end{equation*}%
has infinitely many solutions; moreover we cannot define infinitesimal
numbers such as%
\begin{equation*}
\frac{1}{\infty }
\end{equation*}%
in a consistent way. For this reason the Euler's "$\infty "$ did not even
got the dignity of "number".

\subsubsection{The Hume's Principle and the cardinal numbers}

Until the XIX century, the idea of number and "counting system" was rooted
not only on the Euclid's principle, but also on the Hume's Principle:

\begin{center}
\textit{The number of elements in }$F$\textit{\ is equal to the number of
elements in }$G$\textit{\ }

\textit{if there is a one-to-one correspondence between }$F$\textit{\ and }$%
G $\textit{.}
\end{center}

In our formalism

\begin{itemize}
\item \textbf{Hume's principle} - \textit{Given two sets }$F$\textit{\ and }$%
G$, \textit{then }$F\cong G$ \textit{if there is a bijection} 
\begin{equation*}
\phi :F\rightarrow G.
\end{equation*}
\end{itemize}

\bigskip

Euclid's principle and Hume's Principle are satified by $(\mathbf{Fin},%
\mathbb{N},\left\vert \cdot \right\vert )$ but they lead to a contradiction
if our universe contains an infinite set.

Cantor had the great idea to drop Euclid's principle and to use only the
relation "$\cong $" suggested by Hume and intruduced the \textit{cardinal
numbers} counting system which we will denote by%
\begin{equation*}
(\mathbf{U}_{\max },\mathbb{\ }\mathbf{Card},\ \left\vert \cdot \right\vert
).
\end{equation*}

This counting system is much reacher than $(\mathbf{U},\mathbb{N}\cup
\left\{ \infty \right\} ,\mathfrak{n})$ since for every set $A$,%
\begin{equation*}
\left\vert \wp \left( A\right) \right\vert >\left\vert A\right\vert
\end{equation*}%
Actually Cantor proved that the set of infinite cardinal number form a
sequence 
\begin{equation*}
\aleph _{0}<\aleph _{1}<...<\aleph _{\beta }<...
\end{equation*}%
where $\beta $ is an ordinal number.

Cardinal numbers do not satisfy Euclid's principle, but retain all the basic
algebraic properties (BAP). Unfortunately however, even if the BAP are
satisfied, the lack of Euclid's Principle gives rise to an algebra very
different from the algebra of $(\mathbf{Fin},\mathbb{N},\left\vert \cdot
\right\vert )$; for example, the equation%
\begin{equation*}
a+x=b\ \ (a\leq b)
\end{equation*}%
does not have a unique solution whenever $a$ is infinite; therefore it is
not always possible to define the difference of two cardinal numbers.
Furthermore, their arithmetic is poor since, given two cardinal numbers $%
\alpha $ and $\beta $, if only one of them is\ infinite, we have that 
\begin{equation*}
\alpha +\beta =\alpha \cdot \beta =\max \left( \alpha ,\beta \right)
\end{equation*}%
Then, also in this case, we cannot define infinitesimal number such as%
\begin{equation*}
\frac{1}{\aleph _{\beta }}.
\end{equation*}

\subsubsection{The ordinal numbers}

Cantor introduced also the notion of ordinal number. We can define a triple $%
(\mathbf{W},\mathbf{Ord},\mathfrak{ord})$ as follows:

\begin{itemize}
\item $\mathbf{W}$ is the class of well ordered sets.

\item $\mathbf{Ord}$ is the family of ordnal numbers.

\item $\forall A\in \mathbf{W}$, $\mathfrak{ord}(A)\in \mathbf{Ord}$ is the
order type of $A.$
\end{itemize}

Similarly, we can define a order relation on sets in $\mathbf{W}$ as follows:

\begin{itemize}
\item \textbf{\ }$A\preceq B$ if and only there exists an injection $\Phi
:A\rightarrow B$ which preserves the order, namely $\forall a_{1},a_{2}\in
A, $ 
\begin{equation*}
a_{1}<a_{2}\Rightarrow \Phi \left( a\right) <\Phi \left( b\right) .
\end{equation*}
\end{itemize}

If we equip the ordinal numbers with the operations introduced by Cantor, $(%
\mathbf{W},\mathbf{Ord},\mathfrak{ord})$ is not a counting system since
these operations are different from the ones given by definitions \ref{S}
and \ref{P}. However, if we use the \textit{natural operations} introduced
by Hessenberg they form a counting system. We will come back on this point
in sections \ref{vno} and \ref{SPO}.

\subsubsection{The numerosities}

We have seen that it is not possible to have a counting theory which
contains infinite sets in and which at the same time preserves both Euclid's
principle and Hume's principle. However, we can give up Humes's principle
and keep Euclid's principle.

\begin{definition}
A counting system $(\mathbf{U},\mathbf{Num,}\mathfrak{num)}$ which preserves
the Euclid's principle is called numerosity theory.
\end{definition}

Euclid's Principle is not only inherent to our idea of number, but also
implies an important algebraic property: consider the equation%
\begin{equation}
a+x=b,\ a\leq b;  \label{e}
\end{equation}%
by virtue of the Comparison Principle, this equation always admits a
solution; Euclid's principle implies that this solution is unique. This fact
allows us to define the class of "signed" numbers $\mathcal{Z}$ which is the
analogue of the set of integers $\mathbb{Z}$. $\mathcal{Z}$ can be
(informally) defined in the following way%
\begin{equation*}
\mathcal{Z}=\mathbf{Num}\cup \left\{ -x\ |\ x\in \mathbf{Num}\right\}
\end{equation*}%
The uniqueness of the solution of (\ref{e}) and the basic algebraic
properties allow to prove that $\mathcal{Z}$ is an ordered integrity domain
and therefore it is possible to define the related field of quotients $%
\mathcal{Q}$, i.e. the set of numbers of the form%
\begin{equation*}
x=\pm \frac{\mathfrak{num}(A)}{\mathfrak{num}(B)};\ B\neq \varnothing
\end{equation*}

$\mathcal{Q}$ turns out to be a non-Archimedean field which derives, like
the field of rational numbers, from the idea of "number of elements" of a
set. A further step leads to a field which contains the real number actually
to a field isomprphic to a field of hyperreal numbers and to the field of
surreal numbers. This point will be seen and discussed in section \ref{CONT}.

The numerosity counting system will be denoted by $(\Lambda _{\flat },%
\mathbf{Num,}\mathfrak{num}),$ where 
\begin{equation}
\Lambda :=\{E\in V_{\omega }(\mathbf{Ato})\ |\ \left\vert E\right\vert <%
\mathbf{\kappa }\};\ \Lambda _{\flat }=\Lambda \backslash \mathbf{Ato}
\label{ran}
\end{equation}%
namely, the sets in $\Lambda _{\flat }$ have accessible cardinality and
finite rank. The latter limitation is necessary as the following proposition
shows:

\begin{proposition}
If $(\mathbf{U},\mathbf{Num,}\mathfrak{num)}$ is a numerosity counting
system, and $V_{\beta }(E)\subset \mathbf{U}\cup \mathbf{Ato}$, then $\beta
=\omega $.
\end{proposition}

\textbf{Proof}: We argue indirectly and we assume that $\mathbf{U}$ contais
a set of infinite rank such as%
\begin{equation*}
A=\{a,(a,a),(a,a,a),...\},\ a\in E
\end{equation*}%
then, taking $B=\{a\}$ we have that%
\begin{equation*}
A\times B=\{(a,a),(a,a,a),(a,a,a,a)...\}\subset A
\end{equation*}%
This fact contrdict the definition of counting system since, by the Euclid's
principle we have that 
\begin{equation*}
\mathfrak{num(}A\times B)<\mathfrak{num(}A)
\end{equation*}%
while by Def. \ref{C}-(\ref{C4},\ref{C5}) we have that%
\begin{equation*}
\mathfrak{num(}A\times B)=\mathfrak{num(}A)\times \mathfrak{num(}B)=%
\mathfrak{num(}A)\times 1=\mathfrak{num(}A)
\end{equation*}%
Hence, we must have $\beta \leq \omega .$ Moreover, since%
\begin{equation*}
rank(A\times B)>\max \{rank(A),rank(B)\},
\end{equation*}
we must have $\beta =\omega .$

$\square $

\begin{remark}
If we weaken the request (\ref{C4}) of Def. \ref{C}, it is possible to
develop a numerosity theory also for sets of infinite rank (see e.g. \cite%
{BDNF1}). However, in this paper, we prefer to avoid the thechnicalites
conned to this choice.
\end{remark}

The existence of a numerosity counting system, namely the consitency of the
Euclid's principle with the notion of counting system, will be proved in
section \ref{CN}. In the next two sections will dig into the very rich
consequences that the Euclid's principle implies.

\section{Numerosities and transfinite numbers\label{NT}}

The first peculiarity of $\mathbf{Num}$ is that this set contains in a
natural way other sets of numbers such as $\mathbf{Ord}$ and $\mathbf{Card}$%
\textbf{. }

\subsection{Numerosities and ordinal numbers\label{vno}}

In this section, we will identify a subset of the numerosities with the
initial segment of cardinality \textbf{$\kappa $} of the class of ordinal
numbers. Let us see how. We remember that by the definition of number, 
\begin{equation}
\mathbf{Num}\subset \mathbf{Ato};  \label{na}
\end{equation}%
by this assumption it makes sense to talk of the numerosity of a set of
numerosities.

\begin{definition}
\label{nuovi ord}The set of the ordinal numerosities (which we will denote
by $\mathbf{Ord}$) is defined as follows: $\beta \in \mathbf{Ord}$ if and
only if%
\begin{equation*}
\beta =\mathfrak{num}\left( \mathbf{O(}\beta \mathbf{)}\right)
\end{equation*}%
where $\forall \beta \in \mathbf{Num}$, $\mathbf{O(}\beta \mathbf{)}%
:=\left\{ x\in \mathbf{O\mathbf{rd}\ |\ }x<\beta \right\} .$
\end{definition}

It is easy to see that $\mathbf{Ord\neq \varnothing }$ since $0=\mathfrak{num%
}\left( \mathbf{\varnothing }\right) =\mathfrak{num}\left( \mathbf{%
O(\varnothing )}\right) \in \mathbf{Ord;}$ moreover

\begin{itemize}
\item if $\beta \in \mathbf{Ord,}$ then $\beta +1=\mathfrak{num}\left( 
\mathbf{O(}\beta \mathbf{)}\cup \left\{ \beta \right\} \right) \in \mathbf{%
Ord}$

\item if $\beta \in \mathbf{Ord,}$ then $\beta :=\mathfrak{num}\left(
\dbigcup\limits_{\gamma <\beta }\mathbf{O(}\gamma \mathbf{)}\right) \in 
\mathbf{Ord}$
\end{itemize}

This construction of the ordinal numerosities is similar to the construction
of Von Neumann ordinals. While for Von Neumann an ordinal number $\beta _{%
\text{\textsc{vn}}}$ is the set of all the ordinal numbers contained in $%
\beta _{\text{\textsc{vn}}}$, an ordinal numerosities $\beta $ is the
numerosity of the set of ordinal numerosities smaller than $\beta .$

Obviously, not all numerosities are ordinal: for example, $\mathfrak{num}%
\left( \mathbb{N}^{+}\right) $ is not ordinal. In fact, if 
\begin{equation}
\mathbf{\alpha }:=\mathfrak{num}\left( \mathbb{N}^{+}\right)  \label{alfa}
\end{equation}%
were an ordinal then: 
\begin{eqnarray*}
\mathbf{\alpha } &=&\mathfrak{num(}\left\{ x\in \mathbf{Ord}\ |\ x<\mathfrak{%
num}\left( \mathbb{N}^{+}\right) \right\} )=\mathfrak{num}(\mathbb{N}) \\
&=&\mathfrak{num}(\mathbb{N}^{+}\cup \left\{ 0\right\} )=\mathbf{\alpha }+1.
\end{eqnarray*}

From now on, we will identify the ordinal numbers with the ordinal
numerosities.

\subsubsection{Operations with ordinal numerosities\label{SPO}}

In this section we will compare the operations between numerosities with the
Cantorian operations between ordinals. Since we use the ordinary symbols $+$
and $\cdot $ for the operations on numerosities, the Cantorian
multiplication and addition\emph{\ }on $\mathbf{Ord}$ will be denoted by $%
\oplus $ and $\odot $. Moreover, we will denote by $\beta ^{\left\langle
\gamma \right\rangle }$ the ordinal exponentiation. We recall that $\beta
^{\left\langle \gamma \right\rangle }$ is defined by induction as follows:

\begin{itemize}
\item (i) $\beta ^{\left\langle 0\right\rangle }=1$

\item (ii) $\beta ^{\left\langle \gamma +1\right\rangle }=\beta
^{\left\langle \gamma \right\rangle }\cdot \beta $

\item (iii) $\beta ^{\left\langle \gamma \right\rangle }=\ \underset{\mathbf{%
Ord}}{\sup }\left\{ \beta ^{\left\langle x\right\rangle }\ |\ x\in \mathbf{O(%
}\gamma \mathbf{)}\right\} $ when $\gamma $ is a limit ordinal\footnote{%
Here, we have used the obvious notation 
\begin{equation*}
\underset{\mathbf{Ord}}{\sup }A:=\min \left\{ \gamma \in \mathbf{Ord}\ |\
\forall x\in A,\ \gamma \geq x\right\}
\end{equation*}%
}.
\end{itemize}

From (i) and (ii), it follows that $\forall n\in \mathbb{N}$, 
\begin{equation}
\omega ^{\left\langle n\right\rangle }=\omega ^{n}.  \label{lina}
\end{equation}
However, if $\gamma \notin \mathbb{N},$ the exponential numerosity $\beta
^{\gamma }$ will be defined in a different way (see section \ref{EN}) and
this fact legitimize the choice of the symbol $\beta ^{\left\langle \gamma
\right\rangle }.$

We recall that each ordinal $\gamma $ can be written in the Cantor normal
form, namely 
\begin{equation*}
\beta =\left( \omega ^{\left\langle j_{n}\right\rangle }\odot b_{n}\right)
\oplus \left( \omega ^{\left\langle j_{n-1}\right\rangle }\odot
b_{n-1}\right) \oplus ...\oplus \left( \omega ^{\left\langle \gamma
_{0}\right\rangle }\odot b_{0}\right) :=\dbigoplus\limits_{k=0}^{n}\left(
\omega ^{\left\langle j_{k}\right\rangle }\odot b_{k}\right)
\end{equation*}%
where $b_{k}\in \mathbb{N}$ and $k_{1}>k_{2}\Rightarrow j_{k_{1}}>j_{k_{2}}$.

Using the Cantor normal form, the natural (or Hessenberg) operations "$+_{%
\text{\textsc{h}}}$" and "$\cdot _{\text{\textsc{h}}}$" are defined as
follows: given%
\begin{equation}
\beta =\dbigoplus\limits_{k=0}^{n}\left( \omega ^{\left\langle
j_{k}\right\rangle }\odot b_{k}\right) ,\ \gamma
=\dbigoplus\limits_{k=0}^{n}\left( \omega ^{\left\langle j_{k}\right\rangle
}\odot c_{k}\right)  \label{bg}
\end{equation}%
(where some coefficient can be null in order to have the same set of $j_{k}$%
's), we have 
\begin{equation*}
\beta +_{\text{\textsc{h}}}\gamma :=\dbigoplus\limits_{k=0}^{n}\left[ \omega
^{\left\langle j_{k}\right\rangle }\odot \left( b_{k}+c_{k}\right) \right]
;\ \beta \cdot _{\text{\textsc{h}}}\gamma :=\dbigoplus\limits_{k,l=0}^{n}%
\left[ \omega ^{\left\langle j_{k}+_{\text{\textsc{h}}}j_{l}\right\rangle
}\odot \left( b_{k}c_{l}\right) \right]
\end{equation*}%
It is interesting and somewhat surprising that the natural operations
coincide with the numerosity operations, namely 
\begin{equation}
\beta +_{\text{\textsc{h}}}\gamma =\beta +\gamma \ \ and\ \ \beta \cdot _{%
\text{\textsc{h}}}\gamma =\beta \cdot \gamma  \label{hess}
\end{equation}

Let us prove this fact.

\begin{lemma}
\label{TT}If $\beta =\dbigoplus\limits_{k=0}^{n}\left( \omega ^{\left\langle
j_{k}\right\rangle }\odot b_{k}\right) ,$ then $\beta
=\dsum\limits_{k=0}^{n}b_{k}\omega ^{\left\langle j_{k}\right\rangle }.$
\end{lemma}

\textbf{Proof}: We set $S_{m}=\left\{ \xi \ |\ \xi
<\dbigoplus\limits_{k=0}^{m}\left( \omega ^{\left\langle j_{k}\right\rangle
}\odot b_{k}\right) \right\} ;$ then, we have that 
\begin{equation*}
\mathbf{O(}\beta \mathbf{)}=\left( S_{n}\backslash S_{n-1}\right) \cup
\left( S_{n-1}\backslash S_{n-2}\right) \cup ....\cup \left( S_{1}\backslash
S_{0}\right) \cup S_{0}
\end{equation*}%
and%
\begin{eqnarray*}
S_{m}-S_{m-1} &=&\left\{ \xi \ |\ \dbigoplus\limits_{k=0}^{m-1}\left( \omega
^{\left\langle j_{k}\right\rangle }\odot b_{k}\right) \leq \xi
<\dbigoplus\limits_{k=m-1}^{m}\left( \omega ^{\left\langle
j_{k}\right\rangle }\odot b_{k}\right) \right\} \\
&=&\left\{ \xi \ |\ \xi <\omega ^{\left\langle j_{m}\right\rangle }\odot
b_{m}\right\}
\end{eqnarray*}%
Therefore%
\begin{eqnarray*}
\mathfrak{num}(S_{m}\backslash S_{m-1}) &=&b_{m}\cdot \mathfrak{num}(\left\{
\xi \ |\ \xi <\omega ^{\left\langle j_{m}\right\rangle }\right\} ) \\
&=&b_{m}\cdot \mathfrak{num}\left( \mathbf{Ord}\left( \omega ^{\left\langle
j_{m}\right\rangle }\right) \right) =\omega ^{\left\langle
j_{m}\right\rangle }
\end{eqnarray*}%
and hence 
\begin{eqnarray*}
\beta &=&\mathfrak{num(}\mathbf{O(}\beta \mathbf{)})=\mathfrak{num}\left(
S_{n}\backslash S_{n-1}\right) +....+\mathfrak{num}\left( S_{1}\backslash
S_{0}\right) +\mathfrak{num}S_{0} \\
&=&b_{n}\omega ^{\left\langle j_{n}\right\rangle }+....+b_{1}\omega
^{\left\langle j_{1}\right\rangle }+b_{0}\omega ^{\left\langle
j_{0}\right\rangle }
\end{eqnarray*}

$\square $

\begin{theorem}
The identities (\ref{hess}) are satisfied.
\end{theorem}

\textbf{Proof}: By Th. \ref{TT}, 
\begin{eqnarray*}
\beta +_{\text{\textsc{h}}}\gamma &=&\dbigoplus\limits_{k=0}^{n}\left[
\omega ^{\left\langle j_{k}\right\rangle }\odot \left( b_{k}+c_{k}\right) %
\right] =\dsum\limits_{k=0}^{n}\left( b_{k}+c_{k}\right) \omega
^{\left\langle j_{k}\right\rangle } \\
&=&\dsum\limits_{k=0}^{n}b_{k}\omega ^{\left\langle j_{k}\right\rangle
}+\dsum\limits_{k=0}^{n}c_{k}\omega ^{\left\langle j_{k}\right\rangle
}=\beta +\gamma
\end{eqnarray*}

Moreover, 
\begin{eqnarray*}
\beta \cdot _{\text{\textsc{h}}}\gamma &=&\dbigoplus\limits_{k,l=0}^{n}\left[
\omega ^{\left\langle j_{k}+_{\text{\textsc{h}}}j_{l}\right\rangle }\odot
\left( b_{k}c_{l}\right) \right] =\dbigoplus\limits_{k=0}^{n}\omega
^{\left\langle j_{k}\right\rangle }\odot
b_{k}+\dbigoplus\limits_{l=0}^{n}\omega ^{\left\langle j_{l}\right\rangle
}\odot c_{l} \\
&=&\dsum\limits_{k=0}^{n}b_{k}\omega ^{\left\langle j_{k}\right\rangle
}+\dsum\limits_{l=0}^{n}c_{l}\omega ^{\left\langle j_{l}\right\rangle
}=\beta \gamma .
\end{eqnarray*}

$\square $

\begin{remark}
Within the theory of ordinal numbers the description of $\beta +_{\text{%
\textsc{h}}}\gamma \ $and$\ \ \beta \cdot _{\text{\textsc{h}}}\gamma $ in
terms of well ordered set is rather involved; particularly the description
of the set whose order type is $\beta \cdot _{\text{\textsc{h}}}\gamma .$ On
the contrary, using ordinal numerosities $\beta \cdot _{\text{\textsc{h}}%
}\gamma $ is easily described as the numerosity of the set $\mathbf{O(}\beta 
\mathbf{)}\times \mathbf{O(}\gamma \mathbf{)}.$
\end{remark}

\subsection{Numerosities and cardinal numbers\label{NC}}

As usual, a cardinal number $\chi $ can be identified with the ordinal number%
\begin{equation*}
\chi _{0}=\min \left\{ \gamma \in \mathbf{Ord}\ |\ \left\vert \mathbf{O(}%
\gamma \mathbf{)}\right\vert \geq \chi \right\}
\end{equation*}%
Thanks to this identification, from now on we will assume that 
\begin{equation}
\mathbb{N}\subset \mathbf{Card}\subset \mathbf{Ord}\subset \mathbf{Num.}
\label{CON}
\end{equation}

If a numerosity coincide with a cardinal number will be called \textit{%
cardinal numerosity}. For example, we have that%
\begin{equation*}
\omega =\aleph _{0}=\mathfrak{num}(\mathbb{N}).
\end{equation*}%
and in general, $\forall j<\kappa ,\ \omega _{j}=\aleph _{j}.$ Also, $%
\forall j\in \mathbf{Ord}$, the numerosities $\beth _{j}$\footnote{%
We recall that the Beth numbers are defined by transfinite recursion as
follows
\par
\begin{itemize}
\item $\beth _{0}:=\aleph _{0}$%
\par
\item $\beth _{\beta +1}:=2^{\beth _{\beta }}$%
\par
\item $\beth _{\beta }:=\sup \left\{ \beth _{\gamma }\ |\ \gamma <\beta
\right\} $ if $\beta $ is a limit number.
\end{itemize}
\par
{}} are well defined. If we assume the Continuum Hypotesis, then, $\beth
_{j}=\omega _{j},$ but this assumption is not relevant for the numerosity
theory. In section \ref{bb}, we will analize the numerosity%
\begin{equation}
\beth _{1}=\mathfrak{num}(\wp (\mathbb{N})).  \label{simon}
\end{equation}

\textbf{Caveat! -} From now on, the symbols $\beth _{j}$'s will denote
numerosities and the operations will dentote the operations in the framework
on numerosities.

\subsection{Three different ways of counting\label{twc}}

In this section, we will analyze the operation of counting from a more
intuitive perspective. In section \ref{lu}, we will formalize this operation
introducing the notion of $\Lambda $-limit which, among the other things,
will allow us to move from discrete to continuous, namely from counting to
measuring.

In everyday life, there are several possible ways of counting the number of
elements of finite sets which, of course, yield the same result. However,
when these ways of counting are formalized and extended to infinite sets,
they may give different counting systems. Basically, there are three
different approaches.

\begin{itemize}
\item The first way of counting consists in associating to each element of a
set an element of another one. If in this way one gets a 1-1 correspondence
and claims that the two sets have the same number of elements. This
intuition corresponds to the equipotency relation and to the Cantorian
theory of cardinal numbers.

\item In the second way of counting, one arranges the elements of a given
set in a row and compares such a row with the sequence of natural numbers.
This intuition leads to the notion of order type and to the theory of
ordinal numbers.

\item However, there exists a third way of counting which consists in
arranging the elements of a given sets into smaller groups to be counted
separately. As we will see in the next section this intuition is strongly
related to the notion of numerosity.
\end{itemize}

Please, note that the three ways of counting discussed above imply more and
more complex logical operations.

\begin{itemize}
\item The first way corresponds to the concept of number of a two years old
kid, who associate numbers to sets of fingers of his hands; \emph{e.g.}, the
number 3 corresponds to the set 
\begin{equation*}
\{\text{index finger},\text{middle finger},\text{ring finger}\}.
\end{equation*}%
\smallskip

\item The second way of counting corresponds to the concept of number of a
four years old child: she/he has already memorized the sequence of the first
natural numbers and she/he is able to count objects by arranging them in a
row.

\smallskip

\item The third way of counting is much more sophisticated and requires
several operations, such as collecting similar objects together, and
comparing different groups. This is the way of counting of a grown child.
\end{itemize}

Clearly, the third way of counting is only possible if the objects of a
given set have a "some feature" that allows us to bunch "similar objects".
So we are lead to a structure formalized by the notion of \textit{label}
that will be considered in the next section.{}

\subsection{The label-lattice\label{lu}}

If we want to formalize the third way of counting to any set, we need a
criterion to collect groups of elements. Informally, we may collect elements
sharing the same "label". Now, let us formalize the notion of label.

\begin{definition}
\label{LU}A \textbf{labelling }is a family of sets $\mathfrak{L}$ which
satisfies the following relations:

\begin{enumerate}
\item \label{a}$\mathfrak{L}\subset \mathbf{Fin;}$

\item \label{b}$\lambda ,\mu \in \mathfrak{L}\Rightarrow \lambda \cap \mu
\in \mathfrak{L}\ \ $and$\ \ \exists \sigma \in \mathfrak{L,\ }\lambda \cup
\mu \subseteq \sigma .$

\item \label{c}$\forall a\in \Lambda $, $\exists \lambda \in \mathfrak{L,\ }%
a\in \lambda .$
\end{enumerate}
\end{definition}

Given a labelling $\mathfrak{L}$, the label of an element $a\in \Lambda $ is
defined as follows:%
\begin{equation}
\ell (a)=\dbigcap \left\{ \mu \in \mathfrak{L}\ |\ a\in \mu \right\} .
\label{cica}
\end{equation}%
The set $\mathfrak{L}=\wp _{\omega }(\Lambda ),\ $is the the maximum
labelling. Actually, there exists infinitely many labellings; now, we will
consider a generic labelling $\mathfrak{L}$ since it is sufficient to our
purposes. In section \ref{SP} we will consider peculiar labellings which
provide the numerosity theory with "special properties".

By Def. \ref{LU}-(\ref{b}), $(\mathfrak{L}\mathfrak{,\mathfrak{\subseteq })}$
can be equipped with a lattice structure by setting%
\begin{equation*}
\lambda \wedge \mu :=\lambda \cap \mu ;\ \lambda \vee \mu :=\dbigcap \left\{
\sigma \in \mathfrak{L}\ |\ \lambda \cup \mu \subseteq \sigma \right\} ;
\end{equation*}%
it we will be called \textbf{label-lattice}. In particular, it is a directed
set; then function $\varphi :\mathfrak{L}\rightarrow R$ is a \textit{net }%
(with values in $R$); the set of such nets will be denoted by $\mathfrak{F}%
\left( \mathfrak{L},R\right) .$ If $R$ is a commutative ring, then also $%
\mathfrak{F}\left( \mathfrak{L},R\right) $ is a commutative ring with the
operations defined by%
\begin{equation*}
\left( \varphi +\psi \right) (\lambda )=\varphi (\lambda )+\psi (\lambda );\
\left( \varphi \cdot \psi \right) (\lambda )=\varphi (\lambda )\cdot \psi
(\lambda ).
\end{equation*}%
Given $A\in \Lambda ,$ the \textbf{counting net }$\varphi _{A}\in \mathfrak{F%
}\left( \mathfrak{L},\mathbb{Z}\right) $\textbf{\ }is defined as follows:%
\begin{equation*}
\varphi _{A}(\lambda )=\left\vert \{x\in A\ |\ x\in \lambda \}\right\vert
=\left\vert A\cap \lambda \right\vert .
\end{equation*}%
Now we extend the subtraction in $\mathbf{Num}$ when $\mathfrak{num}(A)<%
\mathfrak{num}(B)$ by setting, 
\begin{equation*}
\mathfrak{num}(A)-\mathfrak{num}(B)=-\mathfrak{num}(B\backslash A^{\prime
}),\ \text{with\ }A^{\prime }\subset A,\ \mathfrak{num}\left( A^{\prime
}\right) =\mathfrak{num}\left( A\right) .
\end{equation*}%
We will denote by $\mathcal{Z}$ the relative ring and we will call it ring
of the \textbf{signed numerosities}.

If $\mathfrak{F}_{\mathcal{Z}}\left( \mathfrak{L},\mathbb{Z}\right) \ $is
the ring generated by the counting nets in $\mathfrak{F}\left( \mathfrak{L},%
\mathbb{Z}\right) ,$ we denote by 
\begin{equation}
J_{\mathcal{Z}}:\mathfrak{F}_{\mathcal{Z}}\left( \mathfrak{L},\mathbb{Z}%
\right) \rightarrow \mathcal{Z}  \label{JZ}
\end{equation}%
the ring homomorphism such that $J_{\mathcal{Z}}(\varphi _{A})=\mathfrak{num}%
(A).$ A number $J_{\mathcal{Z}}(\varphi )\in \mathcal{Z}$ can be seen as a
sort of limit of the net $\lambda \mapsto \varphi (\lambda );$ hence, it is
natural to employ the following notation\footnote{%
This kind of limit generalizes the $\alpha $-limit defined in \cite{BDN2018}
and, as we will see, agrees with the notion of $\Lambda $-limit used in
other papers such as \cite{BLL}.}%
\begin{equation*}
\lim_{\lambda \uparrow \Lambda }\varphi (\lambda )=J_{\mathcal{Z}}(\varphi ).
\end{equation*}%
In order to distinguish the above limit (which we will call $\Lambda $-%
\textit{limit}) from the Cauchy$\ $limit of a net, we have used the notation
"$\lambda \uparrow \Lambda $" rather then "$\lambda \rightarrow \Lambda $".
Hence the numerosity of a set can be expressed as follows: 
\begin{equation}
\mathfrak{num}(A)=\lim_{\lambda \uparrow \Lambda }\ \left\vert A\cap \lambda
\right\vert .  \label{lim}
\end{equation}%
This notation is very significant since $\mathfrak{num}(A)$ appears as the
limit of the cardinality of finite subsets of $A$ that grow with the growth
of $\lambda .$ Furthermore, the idea of $\Lambda $-limit expresses well the
intuitive idea exposed in section \ref{twc}: in order to count the elements
of a large set, it is convenient to count the elements of small groups and
unify the result. Finally, since $J_{\mathbb{Z}}$ is a ring homomorphism,
the $\Lambda $-limit satisfies the some of the properties of the Cauchy$\ $%
limit:

\begin{itemize}
\item if there exists $\lambda _{0}$ such that $\forall \lambda \geq \lambda
_{0},$ $\varphi (\lambda )=\psi (\lambda ),$ then%
\begin{equation*}
\lim_{\lambda \uparrow \Lambda }\varphi (\lambda )=\lim_{\lambda \uparrow
\Lambda }\psi (\lambda )
\end{equation*}

\item if $C_{z}(\lambda )$ is the net identically equal to $z\in \mathbb{Z}$%
, then 
\begin{equation*}
\lim_{\lambda \uparrow \Lambda }C_{z}(\lambda )=z
\end{equation*}

\item for every $\varphi ,\psi \in \mathfrak{F}_{\text{\textsc{c}}}\left( 
\mathfrak{L},\mathbb{Z}\right) ,$ 
\begin{eqnarray*}
\lim_{\lambda \uparrow \Lambda }\left[ \varphi (\lambda )\pm \psi (\lambda )%
\right] &=&\lim_{\lambda \uparrow \Lambda }\varphi (\lambda )\pm
\lim_{\lambda \uparrow \Lambda }\psi (\lambda ) \\
\lim_{\lambda \uparrow \Lambda }\left[ \varphi (\lambda )\cdot \psi (\lambda
)\right] &=&\lim_{\lambda \uparrow \Lambda }\varphi (\lambda )\cdot
\lim_{\lambda \uparrow \Lambda }\psi (\lambda )
\end{eqnarray*}

\item for every $\varphi ,\psi \in \mathfrak{F}_{\text{\textsc{c}}}\left( 
\mathfrak{L},\mathbb{Z}\right) ,$ 
\begin{equation*}
\varphi (\lambda )\geq \psi (\lambda )\Rightarrow \lim_{\lambda \uparrow
\Lambda }\varphi (\lambda )\geq \lim_{\lambda \uparrow \Lambda }\psi
(\lambda ).
\end{equation*}
\end{itemize}

At this point, the notion of $\Lambda $-limit might appear technically
irrelevant, but it will play an important role when it will be extended to
the ring of $\mathbb{R}$-valued nets (see Sec. \ref{NNSA}).

\subsection{The Hume principle revisited}

The notion of $\Lambda $-limit suggests the following definition:

\begin{definition}
\label{CC}If $A\in \Lambda $ and $B\subseteq \Lambda ,$ a bijective map 
\begin{equation*}
\Phi :A\rightarrow B
\end{equation*}%
is called \textbf{comparison map} if 
\begin{equation*}
\lim_{\lambda \uparrow \Lambda }\ \left\vert A\cap \lambda \right\vert
=\lim_{\lambda \uparrow \Lambda }\ \left\vert \Phi (A)\cap \lambda
\right\vert .
\end{equation*}
\end{definition}

From this definition, immediately follow that for every comparison map $\Phi
:A\rightarrow \Lambda $ 
\begin{equation*}
\mathfrak{num}(\Phi (A))=\mathfrak{num}(A).
\end{equation*}%
Then the following facts hold:

\begin{proposition}
\label{20}Let $\Phi :A\rightarrow B$ be a map which eventually preservrs the
labels, namely $\forall a\in A$ 
\begin{equation}
\ell (\Phi (a))=\ell (a)
\end{equation}%
then $\Phi $ is a comparison map and $\mathfrak{num}(\Phi (A))=\mathfrak{num}%
(A)$
\end{proposition}

\textbf{Proof}: Trivial.

$\square $

\begin{proposition}
\label{strega}Let $\Phi :A\rightarrow B$ be a bijective map such that $%
\forall a\in A,$%
\begin{equation}
\ell (\Phi (a))\cap B=\ell (a)
\end{equation}%
then $\Phi $ is a comparison map and $\mathfrak{num}(\Phi (A))=\mathfrak{num}%
(A).$
\end{proposition}

\textbf{Proof}: We have that%
\begin{eqnarray*}
\left\vert B\cap \lambda \right\vert &=&\left\vert \{x\in A\ |\ \Phi (x)\in
\lambda \}\right\vert =\left\vert \{x\in A\ |\ \ell (\Phi (x))=\lambda
\}\right\vert \\
&=&\left\vert \{x\in A\ |\ \ell (x)=\lambda \}\right\vert =\left\vert A\cap
\lambda \right\vert
\end{eqnarray*}

$\square $

\bigskip

Now let us see the relation of $\mathbf{Num}$ with respect to the Hume
Principle. The Hume Principle (HP) can be applied to the numerosities and to
the ordinal numbers provided that we restrict the class of permitted maps,
namely, we have the following situation:

\begin{itemize}
\item HP for Cardinals: $\left\vert A\right\vert =\left\vert B\right\vert $
if and only if there is a bijection $\Phi :A\rightarrow B.$

\item HP for Ordinals: if $A$ and $B$ are two well ordered sets, then $%
\mathfrak{ord}\left( A\right) =\mathfrak{ord(}B)$ if and only if there is a
bijection $\Phi :A\rightarrow B$ which respects the order, namely $\forall
a\in A,\forall b\in B,$ $a<b\Rightarrow \Phi \left( a\right) <\Phi \left(
b\right) .$

\item HP for Numerosities: if $A,B\in \Lambda ,$ then $\mathfrak{num}\left(
A\right) =\mathfrak{num(}B)$ if there exists comparison map $\Phi
:A\rightarrow B.$
\end{itemize}

\section{Numerosities and the continuum\label{CONT}}

In the previous section we have compared the numerosity with the other
numbers used to "count" infinite sets. In this section we will examine the
relation of the numerosities with the numbers used to "measure" continuous
magnitudes. In particular we will relate the numerosity to a peculiar field
of \textbf{hyperreal numbers} (see e.g. \cite{keisler76}) and to the field
of \textbf{surreal numbers} (see e.g. \cite{Conway72}).

\subsection{Euclidean numbers\label{EuN}}

The Euclidean line is fundamental not only for geometry, but also for
analysis and applied mathematics because, once the origin and the unit
element have been chosen, all magnitudes can be represented by its points
that we will call \textbf{Euclidean numbers}. The Euclidean line is
generally identified with the real line, but this identification seems too
restrictive to us since infinite and infinitesimal magnitudes cannot be
described by real numbers. Actually the existence of infinitesimal numbers
has been one of the main problem in all hystory of matematics. We recall
some recent studies on this argument: \cite{K12,K12+,Bk12,Eh06}. In this
paper, we propose a vision of the Euclidean line strictly related to the
numerosities.

\begin{definition}
\label{EL}The Euclidean line $\mathbb{E}$ is a real closed \footnote{%
A field is called\textit{\ real closed} if every polynomial of odd degree
has at least one root.} field which contains the numerosities and such that $%
\forall \xi \in \mathbb{E}$, $\exists \zeta \in \mathbf{Num}$ such that%
\begin{equation}
\left\vert \xi \right\vert \,\leq \zeta .  \label{toto}
\end{equation}
\end{definition}

In this section we will construct the Euclidean via the $\Lambda $-limit and
we will examine some of its properties. In particular, we will see that it
is unique up to isomorphism.

We recall that recently other paths approaching non Archimedan fields have
been investigated (see e.g. \cite{BK},\ \cite{GK} and their references).

\subsubsection{Numerosities and hyperreal numbers\label{NNSA}}

Probably the most relevant property of every numerosity system is that it is
isomorphic to a peculiar subset of the hypernatural numbers $\mathbb{N}^{%
\mathbb{\circledast }}$ as defined in Nonstandard Analysis (NSA), provided
that the hyperreal field $\mathbb{R}^{\mathbb{\circledast }}$ is chosen in a
suitable way. In this section we will construct $\mathbb{R}^{\mathbb{%
\circledast }}$ exploiting the ring of signed numerosities and the notion of 
$\Lambda $-limit.

First of all we need the following lemma:

\begin{lemma}
\label{kappa}If $\varphi \in \mathfrak{F}_{\mathcal{Z}}\left( \mathfrak{L},%
\mathbb{Z}\right) ,$ then there exists a set $K\in \Lambda $ such that $%
\forall \lambda \in \mathfrak{L}$%
\begin{equation*}
\varphi \left( \lambda \right) =\varphi \left( K\cap \lambda \right)
\end{equation*}
\end{lemma}

\textbf{Proof}: By definition, every $\varphi \in \mathfrak{F}_{\mathcal{Z}%
}\left( \mathfrak{L},\mathbb{Z}\right) $ can be written as follows:%
\begin{equation*}
\varphi \left( \lambda \right) =\dsum\limits_{m=0}^{n}z_{m}\varphi
_{A_{m}}\left( \lambda \right) ,\ \ z_{m}\in \mathbb{Z},\ A_{m}\in \Lambda
\end{equation*}%
Since $\varphi _{A_{m}}$ is a counting function, then $\varphi
_{A_{m}}\left( \lambda \right) =\left\vert A_{m}\cap \lambda \right\vert
=\varphi _{A_{m}}\left( A_{m}\cap \lambda \right) $ and hence 
\begin{equation*}
\varphi \left( \lambda \right) =\varphi \left( K\cap \lambda \right)
\end{equation*}%
where $K=A_{1}\cup ....\cup A_{n}.$

$\square $

We set%
\begin{equation*}
\mathfrak{F}_{\text{\textsc{b}}}\left( \mathfrak{L},\mathbb{R}\right)
=\left\{ \varphi \in \mathfrak{F}\left( \mathfrak{L},\mathbb{R}\right) 
\mathfrak{\ }|\ \exists K\in \Lambda ,\ \varphi \left( K\cap \lambda \right)
=\varphi \left( \lambda \right) \right\} ;
\end{equation*}%
it is immediate to see that $\mathfrak{F}_{\text{\textsc{b}}}\left( 
\mathfrak{L},\mathbb{R}\right) $ is a subring of $\mathfrak{F}\left( 
\mathfrak{L},\mathbb{R}\right) ;$ also $\mathfrak{F}_{\mathcal{Z}}\left( 
\mathfrak{L},\mathbb{Z}\right) $ is a subring of $\mathfrak{F}\left( 
\mathfrak{L},\mathbb{R}\right) ;$ now, we denote by $\mathfrak{F}_{\mathbb{E}%
}\left( \mathfrak{L},\mathbb{R}\right) $ the subring of $\mathfrak{F}\left( 
\mathfrak{L},\mathbb{R}\right) $ generated by $\mathfrak{F}_{\text{\textsc{b}%
}}\left( \mathfrak{L},\mathbb{R}\right) $ and $\mathfrak{F}_{\mathcal{Z}%
}\left( \mathfrak{L},\mathbb{Z}\right) .$

\begin{theorem}
\label{NSA}Given a numerosity theory $(\Lambda ,\ \mathbf{Num,\ }\mathfrak{%
num)}$ and a labelling $\mathfrak{L}$, there is an ordered field $\mathbb{E}%
\subset \mathbf{Ato}$ and a surjective ring homomorphism 
\begin{equation}
J:\mathfrak{F}_{\mathbb{E}}\left( \mathfrak{L},\mathbb{R}\right) \rightarrow 
\mathbb{E}  \label{JE}
\end{equation}%
such that

\begin{itemize}
\item $\forall A\in \Lambda $, 
\begin{equation*}
J\left( \varphi _{A}\right) =\mathfrak{num}(A).
\end{equation*}

\item $\mathbf{Num}\subset \mathbb{E}$ and the operations $+$ and $\cdot $
coincide;

\item $\mathbb{R}\subset \mathbb{E}$ and the operations $+$ and $\cdot $
coincide.
\end{itemize}
\end{theorem}

\textbf{Proof: }Let $J_{\mathcal{Z}}$ be the homomorphism defined by (\ref%
{JZ}). It is easy to see that set $\ker (J_{\mathcal{Z}})$ is a prime ideal
in $\mathfrak{F}_{\mathcal{Z}}\left( \mathfrak{L},\mathbb{Z}\right) ;$ hence
the set 
\begin{equation*}
\mathcal{I}_{\mathbb{E}}=\left\{ \varphi \psi \mathfrak{\ }|\ \varphi \in
\ker (J_{\mathcal{Z}}),\psi \in \mathfrak{F}_{\text{\textsc{b}}}\left( 
\mathfrak{L},\mathbb{R}\right) \right\}
\end{equation*}%
is an ideal in $\mathfrak{F}_{\mathbb{E}}\left( \mathfrak{L},\mathbb{R}%
\right) .$ Actually, $\mathcal{I}_{\mathbb{E}}$ is a maximal ideal in $%
\mathfrak{F}_{\text{\textsc{b}}}\left( \mathfrak{L},\mathbb{R}\right) $
since $\mathfrak{F}_{\text{\textsc{b}}}\left( \mathfrak{L},\mathbb{R}\right)
/\mathcal{I}_{\mathbb{E}}$ is a field. Let us check this fact. Take $\left[
\varphi \right] \in \mathfrak{F}_{\mathbb{E}}\left( \mathfrak{L},\mathbb{R}%
\right) /\mathcal{I}_{\mathbb{E}},\ \left[ \varphi \right] \neq 0;\ $we need
to prove that $\left[ \varphi \right] $ has an inverse; let $\chi \in 
\mathfrak{F}_{\mathcal{Z}}\left( \mathfrak{L},\mathbb{Z}\right) $ be the
characteristic function of $\varphi ^{-1}(0),$ then 
\begin{equation*}
\forall \lambda \in \mathfrak{L,\ }\varphi (\lambda )\cdot \chi (\lambda
)=0\ \text{and}\ \varphi (\lambda )+\chi (\lambda )\neq 0
\end{equation*}%
Since $\ker (J_{\mathcal{Z}})$ is a prime ideal and $\varphi (\lambda )\cdot
\chi (\lambda )=0,\ $then $\chi \in \ker (J_{\mathcal{Z}})\subset \mathcal{I}%
_{\mathbb{E}};$ thus%
\begin{equation*}
\left[ \varphi \right] \cdot \left[ \frac{1}{\varphi +\chi }\right] =\left[
\varphi +\chi \right] \cdot \left[ \frac{1}{\varphi +\chi }\right] =1.
\end{equation*}%
So, $\mathfrak{F}_{\text{\textsc{b}}}\left( \mathfrak{L},\mathbb{R}\right) /%
\mathcal{I}_{\mathbb{E}}$ is a field and the projection%
\begin{equation*}
\Pi :\mathfrak{F}_{\mathbb{E}}\left( \mathfrak{L},\mathbb{R}\right)
\rightarrow \mathfrak{F}_{\mathbb{E}}\left( \mathfrak{L},\mathbb{R}\right) /%
\mathcal{I}_{\mathbb{E}}
\end{equation*}%
is a ring homomorphism. Now, we define a field $\mathbb{E}\subset \mathbf{Ato%
}$ isomorphic to $\mathfrak{F}_{\mathbb{E}}\left( \mathfrak{L},\mathbb{R}%
\right) /\mathcal{I}_{\mathbb{E}}$. In order to do this, it is sufficient to
take an injective map%
\begin{equation}
\Psi :\mathfrak{F}_{\mathbb{E}}\left( \mathfrak{L},\mathbb{R}\right) /%
\mathcal{I}_{\mathbb{E}}\rightarrow \mathbf{Ato}  \label{psi}
\end{equation}%
such that $\forall A\in \Lambda ,$\ $\Psi \left( \varphi _{A}\right) =%
\mathfrak{num}\left( A\right) $ and $\forall r\in \mathbb{R},$ $\Psi \left( %
\left[ C_{r}\right] \right) =r$ (here $C_{r}\in \mathfrak{F}_{\mathbb{E}%
}\left( \mathfrak{L},\mathbb{R}\right) $ is a net identically equal to $r$
and $\beta ).\ \mathbb{E}=\func{Im}(\Psi )$ is naturally equipped with a
field structure by setting $a+b=\Psi \left( \Psi ^{-1}(a)+\Psi
^{-1}(b)\right) $ and $ab=\Psi \left( \Psi ^{-1}(a)\cdot \Psi
^{-1}(b)\right) .$

In conclusion, we have constructed a ring homomorphism%
\begin{equation*}
J:=\Psi \circ \Pi :\mathfrak{F}_{\mathbb{E}}\left( \mathfrak{L},\mathbb{R}%
\right) \rightarrow \mathbb{E}
\end{equation*}%
which satisfies the requests of the theorem.

$\square $

Thanks to the homomorphism (\ref{JE}), the notion of $\Lambda $-limit can be
extended to every net in $\mathfrak{F}_{\mathbb{E}}\left( \mathfrak{L},%
\mathbb{R}\right) $ by setting%
\begin{equation*}
\lim_{\lambda \uparrow \Lambda }\varphi (\lambda )=J_{\mathbb{E}}(\varphi ).
\end{equation*}%
Clearly, it satisfies the properties listed at the end of Sec. \ref{lu} and
the following one which is not shared by the Cauchy limit:

\begin{itemize}
\item if there exists $\lambda _{0}$ such that $\forall \lambda \geq \lambda
_{0},$ $\varphi (\lambda )>0,$ then%
\begin{equation*}
\lim_{\lambda \uparrow \Lambda }\varphi (\lambda )>0.
\end{equation*}
\end{itemize}

By Th. \ref{NSA} and by well known results (see e.g. \cite{BDN2005}), $%
\mathbb{E}$ is a hyperreal field. Then, we can exploit the usual notation
and techniques of Nonstandard Analysis:

\begin{itemize}
\item if $\xi ,\zeta \in \mathbb{E}$ we set $\xi \sim \zeta $ if $\xi -\zeta 
$ is infinitesimal;

\item if $\xi \in \mathbb{E}$ is a bounded number then, the standard part of 
$\xi $, $st(\xi ),$ is the only real number $r\sim \xi $. If $\xi $ is
unbounded, then we write $st(\xi )=\pm \infty .$

\item every real function $f\in \mathfrak{F}\left( \mathbb{R},\mathbb{R}%
\right) $ can be extended to $\mathbb{E}$ by setting, for every $\xi
:=\lim_{\lambda \uparrow \Lambda }\varphi (\lambda )$%
\begin{equation}
f^{\circledast }\left( \xi \right) :=\lim_{\lambda \uparrow \Lambda
}f(\varphi (\lambda )).  \label{amy1}
\end{equation}%
As usual, when the meaning is clear from the contest, we omit the $%
"^{\circledast }"$ and we will simply write $f(\xi ).$

\item if $E_{\lambda }\in $ $V_{n}(\mathbb{R})$ is a net of sets, then their 
$\Lambda $-limit is defined by induction over $n$ as follows: if $n=1,$ and $%
E_{\lambda }\subset \mathbb{R}$ 
\begin{equation*}
\lim_{\lambda \uparrow \Lambda }E_{\lambda }:=\left\{ \lim_{\lambda \uparrow
\Lambda }\varphi (\lambda )\ |\ \forall \lambda ,\ \varphi (\lambda )\in
E_{\lambda }\right\} \in V_{1}(\mathbb{E})
\end{equation*}%
and if $E_{\lambda }\in $ $V_{n}(\mathbb{R})$%
\begin{equation*}
\lim_{\lambda \uparrow \Lambda }E_{\lambda }:=\left\{ \lim_{\lambda \uparrow
\Lambda }\varphi (\lambda )\ |\ \forall \lambda ,\ \varphi (\lambda )\in
E_{\lambda }\right\} \in V_{n+1}(\mathbb{E}).
\end{equation*}

\item the hyperreal triple $(\circledast ,\mathbb{R},\mathbb{R}^{\circledast
})$ induces a nonstandard universe $(\circledast ,V_{\omega }(\mathbb{R}%
),V_{\omega }(\mathbb{E}))$ in the sense of Keisler (see \cite{keisler76});
the map 
\begin{equation}
\circledast :V_{\omega }(\mathbb{R})\rightarrow V_{\omega }(\mathbb{E})
\label{aida}
\end{equation}%
is defined as follows:%
\begin{equation*}
A^{\circledast }:=\left\{ \lim_{\lambda \uparrow \Lambda }\varphi (\lambda
)\ |\ \forall \lambda ,\ \varphi (\lambda )\in A\right\}
\end{equation*}%
Notice that this definition is equivalent to the following:%
\begin{equation*}
A^{\circledast }=\lim_{\lambda \uparrow \Lambda }C_{A}(\lambda )
\end{equation*}%
where $C_{A}(\lambda )$ is the net identically equal to $A.$

\item If $K\in \Lambda ,$ and $\left\{ r_{k}\right\} _{k\in K}$ is a real
net (i.e. $r_{k}\in \mathbb{R}$), then the \textbf{hyperfinite sum }is
defined as follows:%
\begin{equation*}
\sum_{k\in K^{\circledast }}r_{k}=\lim_{\lambda \uparrow \Lambda }\left(
\sum_{k\in K\cap \lambda }r_{k}\right) .
\end{equation*}
\end{itemize}

\subsubsection{The structure of Euclidean line}

If $K\in \mathbf{Inf}$, we set%
\begin{equation}
\mathbb{R}[K]:=\left\{ \lim_{\lambda \uparrow \Lambda }\varphi (\lambda \cap
K\mathbf{)}\ |\ \varphi \in \mathfrak{F}\left( \mathfrak{L},\mathbb{R}%
\right) \right\} =\left\{ \sum_{k\in K^{\circledast }}r_{k}\ |\ r_{k}\in 
\mathfrak{F}(K^{\circledast },\mathbb{R})\right\}  \notag
\end{equation}%
by virtue of our construction, $\mathbb{R}[K]$ is a hyperreal field.

\begin{remark}
Actually every hyperreal field obtained by an ultrapower $\mathbb{R}^{I}/%
\mathcal{U}$ is isomorphic to some $\mathbb{R}[K]$ provided that its
cardinality is less than $\kappa .$ The other hyperreal fields can be
obtained as the inductive limit of a suitable set of $\mathbb{R}[K]^{\prime
} $s.
\end{remark}

\begin{theorem}
The hyperreal field defined in Th. \ref{NSA} is isomorphic the Euclidean
line as defined by Def. \ref{EL}.
\end{theorem}

\textbf{Proof}. Since $\mathbb{E}=\mathbb{R}^{\circledast }$ is a hyperreal
field, $\forall \xi \in \mathbb{E}$, there exists $\nu \in \mathbb{N}%
^{\circledast }$ such that$\ \left\vert \xi \right\vert \leq \nu .$ By lemma %
\ref{kappa}, $\nu \in \mathbb{R}[K]$ for some $K\in \Lambda .$ Hence%
\begin{equation*}
\left\vert \xi \right\vert \leq \nu \leq \mathfrak{num}(K).
\end{equation*}%
Moreover, since $\mathbb{E}$ is a hyperreal field, it is real closed and $%
\kappa $-saturated. Hence $\mathbb{E}$ is a real closed $\kappa $-saturated
field of cardinality $\kappa $; hence by well known results, it is unique up
to isomorphisms.

$\square $

The next theorem characterizes the Euclidean numbers as hyperfinite sums of
real numbers.

\begin{theorem}
\label{s1}For every $\xi \in \mathbb{E}$, there exists a family of real
numbers $\{r_{k}\}_{k\in K}$, $K\in \Lambda $ such that%
\begin{equation*}
\xi =\sum_{k\in K^{\circledast }}r_{k}.
\end{equation*}
\end{theorem}

\textbf{Proof}. Let $\left\{ \lambda _{j}\right\} _{j\in \mathbf{Ord}}$ be a
well ordering of $\mathfrak{L}$. Given $\xi =\lim_{\lambda \uparrow \Lambda
}\varphi _{\xi }(K\cap \lambda )\in \mathbb{E}$, we set%
\begin{equation*}
r_{0}:=\varphi _{\xi }(K\cap \lambda _{0})
\end{equation*}
\begin{equation}
r_{j}=\varphi _{\xi }(K\cap \lambda _{j})-\sum_{k\in \lambda _{j}\cap
K,k<j}r_{k}  \label{s2}
\end{equation}%
Hence, $\forall \lambda \in \mathfrak{L}$,%
\begin{equation*}
\varphi _{\xi }(K\cap \lambda _{j})=\sum_{k\in \lambda _{j}\cap K}r_{k}
\end{equation*}%
The conclusion follows taking the $\Lambda $-limit of both sides.

$\square $

By our construction, we have that%
\begin{equation*}
\mathbb{E}\subseteq \mathbf{Ato}
\end{equation*}%
This fact implies that 
\begin{equation*}
V_{\omega }(\mathbb{E})\subseteq \Lambda .
\end{equation*}%
The map%
\begin{equation*}
\circledast :V_{\omega }(\mathbb{R})\rightarrow V_{\omega }(\mathbb{E})
\end{equation*}%
is a nonstandard Universe in the sense of Keisler (\cite{keisler76}, Def.
15.8). If we assume $\mathbb{E}=\mathbf{Ato,}$ and hence 
\begin{equation}
V_{\omega }(\mathbb{E})=\Lambda ,  \label{ANS}
\end{equation}%
we are in the usual framework of NSA. Nevertheless, we can use nonstandard
methods even if we do not assume (\ref{ANS}).

\begin{theorem}
Given a hyperreal field $\mathbb{R}[K],$ there is $\mathbb{\gamma }\in 
\mathbf{Ord}$ such that%
\begin{equation*}
\mathbb{R}[K]\subseteq \mathbb{R}\left( \mathbb{\gamma }\right) :=\left\{
\lim_{\lambda \uparrow \Lambda }\varphi (\left\vert \lambda \cap \mathbf{O(}%
\gamma \mathbf{)}\right\vert )\ |\ \varphi \in \mathfrak{F}_{\mathbb{E}%
}\left( \mathbb{N},\mathbb{R}\right) \right\} .
\end{equation*}
\end{theorem}

\textbf{Proof}: By the comparison principle there exists $\gamma \in \mathbf{%
Ord}$ and $K^{\prime }\subset \mathbf{O}\left( \gamma \right) $ such that $%
\mathfrak{num}(K^{\prime })=\mathfrak{num}(K).$ Then, denoting by $\chi _{A}$
the characteristic function of $A,$ it holds 
\begin{equation*}
\sum_{k\in K\cap \lambda }\chi _{\{k\}}(k)=\left\vert K\cap \lambda
\right\vert =\left\vert K^{\prime }\cap \lambda \right\vert =\sum_{k\in
K^{\prime }\cap \lambda }\chi _{\{k\}}(k)
\end{equation*}%
and hence%
\begin{equation*}
\sum_{p\in K\cap \lambda }r_{k}=\sum_{k\in K\cap \lambda }r_{k}\chi
_{\{k\}}(k)=\sum_{k\in K^{\prime }\cap \lambda }r_{k}\chi
_{\{k\}}(k)=\sum_{p\in K^{\prime }\cap \lambda }r_{k}
\end{equation*}%
So, for every $\xi \in \mathbb{R}[K]$, we have that 
\begin{equation*}
\xi =\lim_{\lambda \uparrow \Lambda }\sum_{p\in \lambda \cap
K}r_{k}=\lim_{\lambda \uparrow \Lambda }\sum_{p\in K^{\prime }\cap \lambda
}r_{k}\in \mathbb{R}\left( \mathbb{\gamma }\right) .
\end{equation*}

$\square $

Thanks to the above theorem, the field of Euclidean numbers can be filtered
as follows:%
\begin{equation}
\mathbb{E}=\dbigcup\limits_{\mathbb{\gamma <\mathbf{\kappa }}}\mathbb{R}%
\left( \mathbb{\gamma }\right) .  \label{amy}
\end{equation}%
This representation of the Euclidean numbers is quite useful; in fact every
Euclidean number can be seen as the $\Lambda $-limit of a net defined on $%
\wp _{\omega }(\mathbf{O}\left( \gamma \right) );$ moreover, every
numerosity is equal to the numerosity of a subset of $\mathbf{O}\left(
\gamma \right) .$

\subsubsection{Exponentiation of numerosities\label{EN}}

Since we have identified the numerosities with hypernatural numbers, if $f$
is a real function and $\beta $ is a numerosity, then $f(\beta )$ is a
hyperreal number and it could be the numerosity of some set. In particular,
if $\beta $ and $\gamma $ are numerosities, then $\beta ^{\gamma }$ is a
numerosity defined in a different way than the ordinal numerosity $\beta
^{\left\langle \gamma \right\rangle }$. Let us investigate this point.

\begin{proposition}
\label{annamaria}

\begin{enumerate}
\item \label{vii}$\forall E\in \Lambda ,$%
\begin{equation*}
\mathfrak{num}\left( \mathbf{\wp }_{\omega }\left( E\right) \right) =2^{%
\mathfrak{num}\left( E\right) };
\end{equation*}

\item \label{v}if $X$ is a finite set, then $\forall E\in \Lambda ,$%
\begin{equation*}
\mathfrak{num}\left( \mathfrak{F}\left( X,E\right) \right) =\mathfrak{num}%
(X)^{\mathfrak{num}\left( E\right) }=\left\vert X\right\vert ^{\mathfrak{num}%
\left( E\right) };
\end{equation*}

\item \label{viii}if $X,E\in \Lambda \backslash \{\varnothing \}$, we set 
\begin{equation*}
\mathfrak{F}_{fin}\left( X,E\right) :=\left\{ f:D\rightarrow E\mid D\in 
\mathbf{\wp }_{\omega }(X)\right\} ;
\end{equation*}%
then, if $a\in E$, 
\begin{equation*}
\mathfrak{num}\left( \mathfrak{F}_{fin}\left( X,E\backslash \left\{
a\right\} \right) \right) =\mathfrak{num}\left( E\right) ^{\mathfrak{num}%
(X)}.
\end{equation*}
\end{enumerate}
\end{proposition}

\textbf{Proof}: (\ref{vii}) - Let $\lambda \in \mathfrak{L};$ $a\in \mathbf{%
\wp }_{\omega }\left( E\right) \cap \lambda $ if and only if $a\subseteq
E\cap \lambda ;$ then 
\begin{equation*}
\left\vert \mathbf{\wp }_{\omega }\left( E\right) \cap \lambda \right\vert
=2^{\left\vert E\cap \lambda \right\vert }
\end{equation*}%
Hence%
\begin{eqnarray*}
\mathfrak{num}\left( \mathbf{\wp }\left( E\right) \right) &=&\lim_{\lambda
\uparrow \Lambda }\left\vert \mathbf{\wp }\left( E\right) \cap \lambda
\right\vert =\lim_{\lambda \uparrow \Lambda }2^{\left\vert E\cap \lambda
\right\vert } \\
&=&2^{\lim_{\lambda \uparrow \Lambda }\left\vert E\cap \lambda \right\vert
}=2^{\mathfrak{num}\left( E\right) }
\end{eqnarray*}

(\ref{v}) - Take $\lambda $ so large that $X\subseteq \lambda .$ $f\in 
\mathfrak{F}\left( X,E\right) \cap \lambda $ if and only if$\ $ $\func{Im}%
f\subseteq E\cap \lambda ;$ then%
\begin{equation*}
\left\vert \mathfrak{F}\left( X,E\right) \cap \lambda \right\vert
=\left\vert X\right\vert ^{\left\vert E\cap \lambda \right\vert }
\end{equation*}%
Hence%
\begin{eqnarray*}
\mathfrak{num}\left( \mathfrak{F}\left( X,E\right) \right) &=&\lim_{\lambda
\uparrow \Lambda }\left\vert \mathfrak{F}\left( X,E\right) \cap \lambda
\right\vert =\lim_{\lambda \uparrow \Lambda }\left\vert X\right\vert
^{\left\vert E\cap \lambda \right\vert } \\
&=&\left\vert X\right\vert ^{\lim_{\lambda \uparrow \Lambda }\left\vert
E\cap \lambda \right\vert }=\left\vert X\right\vert ^{\mathfrak{num}\left(
E\right) }
\end{eqnarray*}

(\ref{viii}) We set 
\begin{equation*}
\widetilde{f}(x)=%
\begin{cases}
f(x) & \text{if}\ x\in X\cap \lambda ; \\ 
\ \ a & \text{otherwise}%
\end{cases}%
\end{equation*}%
the map $f\in \mathfrak{F}_{fin}\left( X\cap \lambda ,\left( E\cap \lambda
\right) \backslash \{a\}\right) \rightarrow \widetilde{f}\in \mathfrak{F}%
\left( X\cap \lambda ,E\cap \lambda \right) $ is a bijection. Moreover, if $%
\widetilde{f}\in \mathfrak{F}\left( X\cap \lambda ,E\cap \lambda \right) $ 
\begin{equation*}
\left\vert \mathfrak{F}_{fin}\left( X\cap \lambda ,E\backslash \{a\}\right)
\cap \lambda \right\vert =\left\vert \mathfrak{F}\left( X\cap \lambda ,E\cap
\lambda \right) \right\vert =\left\vert X\cap \lambda \right\vert
^{\left\vert E\cap \lambda \right\vert }
\end{equation*}%
Then%
\begin{eqnarray*}
\mathfrak{num}\left( \mathfrak{F}_{fin}\left( X,E\backslash \{a\}\right)
\right) &=&\lim_{\lambda \uparrow \Lambda }\ \left\vert \mathfrak{F}%
_{fin}\left( X,E\backslash \{a\}\right) \cap \lambda \right\vert \\
&=&\lim_{\lambda \uparrow \Lambda }\ \left\vert \mathfrak{F}_{fin}\left(
X\cap \lambda ,\left( E\cap \lambda \right) \backslash \{a\}\right)
\right\vert \\
&=&\lim_{\lambda \uparrow \Lambda }\ \left\vert \mathfrak{F}\left( X\cap
\lambda ,E\cap \lambda \right) \right\vert \\
&=&\lim_{\lambda \uparrow \Lambda }\ \left\vert E\cap \lambda \right\vert
^{\left\vert X\cap \lambda \right\vert }=\mathfrak{num}\left( E\right) ^{%
\mathfrak{num}\left( X\right) }.
\end{eqnarray*}

$\square $

The cardinality of the continuum is $\left\vert \wp (\mathbb{N})\right\vert $
but we have seen that $2^{\mathbf{num}(\mathbb{N})}$ is the numerosity of
the denumerable set $\wp _{\omega }(\mathbb{N}).$ Actually the numerosity of 
$\wp (\mathbb{N})$ is given by $\beth _{1}$ as defined in section \ref{NC}.
Thus, in the theory of numerosities, there are three different kind of
exponentiation: if $\beta =\mathfrak{num}\left( E\right) $ and $\gamma =%
\mathfrak{num}\left( C\right) $, we have:

\begin{itemize}
\item \textbf{hyperreal exponentiation:}%
\begin{equation*}
\beta ^{\gamma }=\mathfrak{num}\left( \mathfrak{F}_{fin}\left( C,B\backslash
\left\{ b\right\} \right) \right) =\lim_{\lambda \uparrow \Lambda }\
\left\vert B\cap \lambda \right\vert ^{\left\vert C\cap \lambda \right\vert
}\ \ (b\in B);
\end{equation*}

\item \textbf{ordinal exponentiation}: 
\begin{equation*}
\beta ^{\left\langle \gamma \right\rangle }=\underset{\mathbf{Ord}}{\sup }%
\left\{ \mathfrak{num}\left( \mathbf{O}\left( \beta ^{\left\langle
x\right\rangle }\right) \right) \ |\ x\in \mathbf{O}\left( \gamma \right)
\right\} ;
\end{equation*}

\item \textbf{cardinal exponentiation}: if $\beta $ and $\gamma $ are
cardinal numerosities,%
\begin{equation*}
\beta ^{\uparrow \gamma }:=\mathfrak{num}\left( \mathfrak{F}\left(
C,B\right) \right) =\lim_{\lambda \uparrow \Lambda }\ \left\vert \mathfrak{F}%
\left( C,B\right) \cap \lambda \right\vert .
\end{equation*}
\end{itemize}

If $\beta $ and $\gamma $ are ordinal numerosities, $\beta ^{\left\langle
\gamma \right\rangle }=\mathfrak{num}\left( \mathbf{Ord}\left( \beta
^{\left\langle x\right\rangle }\right) \right) ;$ if $\beta $ and $\gamma $
are cardinal numerosities, $\beta ^{\uparrow \gamma }=\mathfrak{num}\left(
B^{C}\right) .$ If $\gamma =n\in \mathbb{N},$ then 
\begin{equation*}
\beta ^{n}=\beta ^{\left\langle n\right\rangle }=\beta ^{\uparrow \gamma }=%
\underset{n\ \text{times}}{\underbrace{\beta \cdot ....\cdot \beta }}
\end{equation*}%
but it is not surprising that these operations give different results when
they are generalized, since they correspond to different ways of counting
infinite sets. We have the following result:

\begin{theorem}
\label{CAC}If $\beta $ and $\gamma $ are ordinal numerosities,%
\begin{equation}
\beta ^{\left\langle \gamma \right\rangle }\leq \beta ^{\gamma }.
\label{lana}
\end{equation}
\end{theorem}

\textbf{Proof}: For $\gamma =0,$ the result is obvious. Now let us argue by
induction: we have that%
\begin{equation*}
\beta ^{\left\langle \gamma +1\right\rangle }=\beta ^{\left\langle \gamma
\right\rangle }\cdot \beta \leq \beta ^{\gamma }\cdot \beta =\beta ^{\gamma
+1}
\end{equation*}%
and if $\gamma $ is a limit ordinal, we have that, $\forall x\,<\gamma ,$%
\begin{equation*}
\beta ^{\gamma }=\mathfrak{num}\left( \mathfrak{F}_{fin}\left( \mathbf{O(}%
\gamma \mathbf{)},\mathbf{O(}\beta \mathbf{)}\backslash \left\{ 0\right\}
\right) \right) >\mathfrak{num}\left( \mathfrak{F}_{fin}\left( \mathbf{O}%
\left( x\right) ,\mathbf{O(}\beta \mathbf{)}\backslash \left\{ 0\right\}
\right) \right) =\beta ^{x}
\end{equation*}%
Hence 
\begin{equation*}
\beta ^{\left\langle \gamma \right\rangle }=\ \underset{\mathbf{Ord}}{\sup }%
\left\{ \beta ^{\left\langle x\right\rangle }\ |\ x\in \mathbf{O(}\gamma 
\mathbf{)}\right\} \leq \underset{\mathbf{Ord}}{\sup }\left\{ \beta ^{x}\ |\
x\in \mathbf{O(}\gamma \mathbf{)}\right\} \leq \beta ^{\gamma }.
\end{equation*}

$\square $

$\bigskip $

\textbf{Example:}

\begin{itemize}
\item $2^{\omega }=\mathfrak{num}\left( \wp _{\omega }\left( \mathbb{N}%
\right) \right) ;$

\item $2^{\left\langle \omega \right\rangle }=\underset{\mathbf{Ord}}{\sup }%
\left\{ 2^{x}\ |\ x\in \mathbb{N}\right\} =\omega ;$

\item $2^{\uparrow \omega }=\mathfrak{num}\left( \wp \left( \mathbb{N}%
\right) \right) =\beth _{1}.$
\end{itemize}

\subsection{Numerosities and surreal numbers\label{NS}}

In this section we will see how the Euclidean numbers can be identified with 
$\mathbf{No}$, namely\textbf{\ }the field of surreal numbers, (see \cite%
{Conway72,Conway,gonsh}) or, to be more precise the pseudoclass of surreal
numbers "born" before the day $\mathbf{\kappa }$. Hence $\mathbf{No}$, using
the terminology introduced in section \ref{pn}, is a pseudoclass. It is well
known that every hyperreal field can be embedded in \textbf{No }and that the
Keisler hyperreal field of cardinality $\kappa $ (see \cite{keisler76}) is
isomorphic to the pseudoclass $\mathbf{No;}$ however, if we consider the
hyperreal field of Euclidean numbers $\mathbb{E}$ (see Def. \ref{EL}), there
is a peculiar isomorphism which relates ordinal numbers, numerosities and
surreal numbers.

In this section, we assume the reader to be familiar with the surreal
numbers; however, we will recall some basic feature of $\mathbf{No}$ and we
will fix the notation.

A surreal number can be represented by a sequence of $+$'s and $-$'s. For
example $+--++$ represents a number and we will write%
\begin{equation*}
x=\left( +--++\right)
\end{equation*}%
In a formal way, a surreal number can be identified with a map 
\begin{equation*}
x:\mathbf{O(}\gamma \mathbf{)}\rightarrow \{+,-\};
\end{equation*}%
the ordinal number $\gamma $ is called \textbf{birthday} of $x$ and it will
be denoted by $\mathbf{b}(x).$

For example, the ordinal numbers (which are supposed to be a sub-pseudoclass
of \textbf{No}) can be represents ad follows:

\begin{itemize}
\item $0=()$

\item $1=(+)$

\item $2=(++)$ etc.

\item $\omega =(++.....)$

\item etc.
\end{itemize}

The surreal numbers can be lexicographically ordered with the convention
that "undefined values" are greater than $-$ and less than $+;$ for example%
\begin{equation*}
(-)<(-+)<(\ \ )<(+-)<(+-+)<(+)<(++-)
\end{equation*}%
We now recall some facts relative to the surreal numbers which we will use
later:

\begin{itemize}
\item The sequences of all $+$'s is can be identified with the ordinal
numbers; hence we may assume $\mathbf{Ord}\subset \mathbf{No;}$

\item $\mathbf{No}$ is a field and the operations + and $\cdot $ restricted
to $\mathbf{Ord}$ agree with the Hessenberg's operations.
\end{itemize}

Now, let us recall the "basic" theorem of the theory of surreal numbers. In
order to do this, it is convenient the following notation: let $A,B\subset 
\mathbf{No}$, ($\left\vert A\right\vert ,\left\vert B\right\vert <\mathbf{%
\kappa }$), be two sets such that $\forall a\in A,\forall b\in B,\ a\,<b;\ $%
then we set%
\begin{equation*}
I(A,B)=\left\{ x\in \mathbf{No}\ |\ \forall a\in A,\forall b\in B,\
a\,<x<b\right\}
\end{equation*}%
namely, $I(A,B)$ is an interval between $A$ and $B.$

\begin{theorem}
(\textbf{Conway simplicity theorem}) Let $A,B\subset \mathbf{No}$, ($%
\left\vert A\right\vert ,\left\vert B\right\vert <\mathbf{\kappa }$), be two
sets such that $\forall a\in A,\forall b\in B,\ a\,<b;\ $then there exists a
uniqie$\ c\in I(A,B)\mathbf{\ }$such that%
\begin{equation*}
\forall x\in I(A,B)\backslash \{c\},\ \mathbf{b}(c)<\ \mathbf{b}(x).
\end{equation*}
\end{theorem}

The number $c$, as defined by the above theorem is denoted by:%
\begin{equation}
c=A|B  \label{cAB}
\end{equation}%
Given a number $x\in \mathbf{No},$ the \textbf{Conway canonical form} of $x$
is given by%
\begin{equation*}
x=L(x)|R(x)
\end{equation*}%
where%
\begin{eqnarray*}
L(x) &=&\left\{ t\in \mathbf{No}\ |\ \text{ }t<x,\ \mathbf{b}(t)<\mathbf{b}%
(x)\right\} \ \text{is the set of "left" numbers born before }x,\  \\
R(x) &=&\left\{ t\in \mathbf{No}\ |\ \text{ }t>x,\ \mathbf{b}(t)<\mathbf{b}%
(x)\right\} \ \text{is the set of "right" numbers born before }x.
\end{eqnarray*}

Now we set 
\begin{equation*}
\mathbf{No}\left( \gamma \right) =\left\{ x\in \mathbf{No}\ |\ x<\gamma
\right\}
\end{equation*}%
and we recall an important result of Van den Dries and Ehrlich \cite%
{Ehrlich01}:

\begin{theorem}
If $\gamma $ is an $\varepsilon $-number, i.e. $\gamma =\omega
^{\left\langle \gamma \right\rangle },$ then $\mathbf{No}\left( \gamma
\right) $ is a subfield of $\mathbf{No.}$
\end{theorem}

In order to proceed, we need to analize some features of the ordinal numbers.

\begin{theorem}
If $\gamma $ is an $\varepsilon $-number, i.e. $\gamma =\omega
^{\left\langle \gamma \right\rangle },$ then $\mathbf{No}\left( \gamma
\right) $ is a subfield of $\mathbf{No.}$
\end{theorem}

\begin{definition}
\label{TETA}A ordinal number $\theta \neq 0$ is called \textbf{indecomposable%
}\footnote{%
Sometimes in the literature, the word "indecomposable" is used as synonimous
of "irreducible"; a number $\gamma $ is irreducible if 
\begin{equation*}
\alpha ,\beta \,<\gamma \Rightarrow \alpha +\beta \,<\gamma .
\end{equation*}%
Here, we use it with a different meaning.} if 
\begin{equation*}
\alpha ,\beta ,\gamma \in \mathbf{Ord}\left( \theta \right) \Rightarrow
\alpha +\beta \gamma <\theta
\end{equation*}%
An ordinal $\gamma \neq 0$ is called \textbf{component} of $\beta $ if there
exist $x,y\in \mathbf{Ord}$ such that%
\begin{equation*}
\gamma +x=\beta \ \ \text{or\ \ }x+\gamma y=\beta .
\end{equation*}
\end{definition}

In order to fix the ideas, we recall that the the indecomposibe ordinals
have the following representation:%
\begin{equation}
\theta _{\gamma }=\omega ^{\left\langle \omega \right\rangle ^{\left\langle
\gamma \right\rangle }},\ \gamma \in \mathbf{Ord}.  \label{teta2}
\end{equation}

The set $\mathbf{Ord}(\theta _{j})$ is closed for the operations + and $%
\cdot $ and every ordinal number $\beta \in \mathbf{Ord}(\theta _{j+1})$ can
be written as follows: 
\begin{equation}
\beta =\dsum\limits_{k=0}^{n}b_{k}\theta _{j}^{k},\ \ b_{k}\in \mathbf{Ord}%
(\theta _{j})  \label{teta}
\end{equation}

\begin{remark}
If $\beta >0,$ there is a unique polynomial of degree $>1$ which allows to
represent the ordinal number $\beta $ in the form (\ref{teta}). In the
context of numerosities this representation of an ordinal number is more
convenient that the Cantor normal form, since the operations involved in it
coincide with the operations defined for the numerosities (and hence for the
Euclidean numbers).
\end{remark}

Now, we denote by $\mathbf{S}(\theta _{\gamma })$ the field generated by $%
\mathbf{O}(\theta _{\gamma })\cup \{\theta _{\gamma }\}$, namely smallest
subfield of $\mathbf{No}$ containing $\mathbf{O}(\theta _{\gamma })\cup
\{\theta _{\gamma }\}.$ Every $\sigma \in \mathbf{S}(\theta _{\gamma })$ can
be written as follows:%
\begin{equation}
\sigma =\frac{\dsum\limits_{k=0}^{n}b_{k}\theta _{\gamma }}{%
\dsum\limits_{k=0}^{n}c_{k}\theta _{\gamma }},\ \ b_{k},c_{k}\in \mathbf{S}%
(\theta _{\gamma })\mathbb{\ }  \label{esse}
\end{equation}%
Moreover, for every number $\mathbb{\beth }_{j}$ we define the field%
\begin{equation*}
\mathbf{S}(\mathbb{\beth }_{j}^{-})=\dbigcup\limits_{\mathbb{\theta }%
_{\gamma }<\mathbb{\beth }_{j}}\mathbf{S}(\mathbb{\theta }_{\gamma })
\end{equation*}%
Since the ordinal numbers are surreal numbers, we may assume that%
\begin{equation*}
\mathbf{S}(\mathbb{\beth }_{j}^{-})\subset \mathbf{No}(\mathbb{\beth }_{j+1})
\end{equation*}%
In fact, if $\theta _{\gamma }=\mathbf{b}(\theta _{\gamma })<\mathbb{\beth }%
_{j+1},\ $then $\mathbf{O}(\theta _{\gamma })\cup \{\theta _{\gamma
}\}\subset \mathbf{No}(\mathbb{\beth }_{j+1})$; since $\mathbb{\beth }_{j+1}$
is an $\varepsilon $-number, $\mathbf{No}(\mathbb{\beth }_{j+1})$ is s
field, and hence $\mathbf{S}(\mathbb{\theta }_{\gamma })\subset \mathbf{No}(%
\mathbb{\beth }_{j+1}).$

The ordinal numbers can be identified also with the Euclidean numbers and by
(\ref{amy1}) and (\ref{amy}), we have that $f\left( \mathbb{\theta }_{\gamma
}\right) \in \mathbb{R}(\mathbb{\theta }_{\gamma })$ for every real function 
$f$; hence $\mathbf{S}(\mathbb{\theta }_{\gamma })\subset \mathbb{R}(\mathbb{%
\theta }_{\gamma })$ and, if $\mathbb{\theta }_{\gamma }<\mathbb{\beth }%
_{j+1}$, $\mathbf{S}(\mathbb{\theta }_{\gamma })\subset \mathbb{R}(\mathbb{%
\beth }_{j}).$ In conclusion,%
\begin{equation*}
\mathbf{S}(\mathbb{\beth }_{j+1}^{-})\subset \mathbb{R}(\mathbb{\beth }_{j}).
\end{equation*}

\begin{lemma}
\label{lulu3}For every $\xi \in \mathbb{R}(\mathbb{\beth }_{j})^{+},$ there
exists $\sigma \in \mathbf{S}(\mathbb{\beth }_{j+1}^{-})^{+}$ such that $%
\sigma <\xi .$
\end{lemma}

\textbf{Proof}: If $\xi \in \mathbb{R}(\mathbb{\beth }_{j})$, then $1/\xi
\in \mathbb{R}(\mathbb{\beth }_{j})$ and by (\ref{toto}), $1/\xi <\mathbb{%
\theta }_{\gamma }$ for some $\mathbb{\theta }_{\gamma }<\mathbb{\beth }%
_{j}. $ Then 
\begin{equation*}
\xi <\frac{1}{\mathbb{\theta }_{\gamma }}\in \mathbf{S}(\mathbb{\beth }%
_{j+1}^{-})^{+}
\end{equation*}

$\square $

\begin{corollary}
\label{lulu4}$\mathbf{S}(\mathbb{\beth }_{j+1}^{-})$ is dense in $\mathbb{R}(%
\mathbb{\beth }_{j})$ with respect to the order topology, namely if $\xi \in 
\mathbb{R}(\mathbb{\beth }_{j}),$ then $\forall \varepsilon \in \mathbb{R}(%
\mathbb{\beth }_{j})^{+},\ \exists \sigma ,\tau \in \mathbf{S}(\mathbb{\beth 
}_{j+1}^{-})$, such that $\sigma \leq \xi \leq \tau $ and $\tau -\sigma \leq
\varepsilon .$
\end{corollary}

\textbf{Proof}:\ It follows from lemma \ref{lulu3} and standard arguments.

$\square $

\begin{definition}
A section of $\mathbf{S}(\mathbb{\beth }_{j}^{-})$ is a pair of non empty
sets $(A,B)$ such that

\begin{itemize}
\item $\forall a\in A,\forall b\in B,\ a\,<b.$

\item $A\cup B=\mathbf{S}(\mathbb{\beth }_{j}^{-})$ or $\exists \sigma \in 
\mathbf{S}(\mathbb{\beth }_{j}^{-}),$ $A\cup \{\sigma \}\cup B=\mathbf{S}(%
\mathbb{\beth }_{j}^{-}).$
\end{itemize}
\end{definition}

The set of all the sections of $\mathbf{S}(\mathbb{\beth }_{j}^{-})$ will be
denoted by $\mathbf{Sec}(\mathbb{\beth }_{j}^{-}).$ Every $\xi \in \mathbb{R}%
(\mathbb{\beth }_{j})$ determines a section in $\mathbf{S}(\mathbb{\beth }%
_{j+1}^{-})$ and we will use the following notation:%
\begin{equation}
\mathcal{L}(\xi )=\left\{ t\in \mathbf{S}(\mathbb{\beth }_{j+1}^{-})\ |\ 
\text{ }t<\xi \right\} ;\ \mathcal{R}(\xi )=\left\{ t\in \mathbf{S}(\mathbb{%
\beth }_{j+1}^{-})\ |\ \text{ }t>\xi \right\} .\   \label{IVcorda}
\end{equation}%
Also the converse is true:

\begin{lemma}
\label{aria2}If $(A,B)$ is a section of $\mathbf{S}(\mathbb{\beth }%
_{j+1}^{-}),$ there exist a uniqe number $\xi \in \mathbb{R}(\mathbb{\beth }%
_{j})$ such that $\left( \mathcal{L}(\xi ),\mathcal{R}(\xi )\right) =(A,B).$
\end{lemma}

\textbf{Proof}: Given $(A,B),$ since $\mathbb{E}$ is $\kappa $-sturated,
there exists $\xi _{0}\in \mathbb{E}$ such that 
\begin{equation*}
\forall a\in A,\forall b\in B,\ a<\xi \,_{0}<b.
\end{equation*}%
However, we cannot conclude that $\xi \,_{0}\in \mathbb{R}(\mathbb{\beth }%
_{j}).$

By definition of Euclidean number, we have that $\forall \sigma \in
A,\forall \tau \in B,$%
\begin{equation*}
\xi _{0}=\lim_{\lambda \uparrow \Lambda }\varphi _{\xi _{0}}(\lambda );\
\sigma =\lim_{\lambda \uparrow \Lambda }\varphi _{\sigma }(\lambda );\ \tau
=\lim_{\lambda \uparrow \Lambda }\varphi _{\tau }(\lambda )
\end{equation*}%
Choosing $\varphi _{\sigma }(\lambda ),\varphi _{\xi _{0}}(\lambda ),\varphi
_{\tau }(\lambda )$ properly, we may assume that $\forall \lambda \in 
\mathfrak{L}$, 
\begin{equation*}
\varphi _{\sigma }(\lambda )<\varphi _{\xi _{0}}(\lambda )<\varphi _{\tau
}(\lambda )
\end{equation*}%
then, in paricular, we have that%
\begin{equation*}
\varphi _{\sigma }(\lambda \cap \mathbf{O}(\beth _{j}))<\varphi _{\xi
_{0}}(\lambda \cap \mathbf{O}(\beth _{j}))<\varphi _{\tau }(\lambda \cap 
\mathbf{O}(\beth _{j}))
\end{equation*}%
Since $\sigma ,\tau \in \mathbb{R}(\mathbb{\beth }_{j}),$ taking the $%
\Lambda $-limit, we get%
\begin{equation*}
\sigma <\lim_{\lambda \uparrow \Lambda }\varphi _{\xi _{0}}(\lambda \cap 
\mathbf{O}(\beth _{j}))<\tau
\end{equation*}%
The conclusion follows taking 
\begin{equation*}
\xi =\lim_{\lambda \uparrow \Lambda }\varphi _{\xi _{0}}(\lambda \cap 
\mathbf{O}(\beth _{j}));
\end{equation*}%
in fact, the uniqueness is guaranteed by Corollary \ref{lulu4}.

$\square $

By the Conway Simplicity Theorem, the map

\begin{equation*}
\{\cdot \}|\{\cdot \}:\mathbf{Sec}(\mathbb{\beth }_{j+1})\rightarrow \mathbf{%
No}(\mathbb{\beth }_{j+1}).
\end{equation*}%
is well defined. So, we can define a map%
\begin{equation}
i:\mathbb{R}(\mathbb{\beth }_{j})\rightarrow \mathbf{No}(\mathbb{\beth }%
_{j+1}),\ \ i(\xi ):=\mathcal{L}(\xi )|\mathcal{R}(\xi )  \label{lara}
\end{equation}

\begin{lemma}
The map $i:\mathbb{R}(\mathbb{\beth }_{j})\rightarrow \mathbf{No}(\mathbb{%
\beth }_{j+1})$ is surjective.
\end{lemma}

\textbf{Proof}: Given $x\in \mathbf{No}(\mathbb{\beth }_{j+1}),$ we set 
\begin{equation*}
\mathcal{L}(x):=\left\{ t\in \mathbf{S}(\mathbb{\beth }_{j+1}^{-})\ |\ t<\xi
\right\} ;\ \mathcal{R}(x)=\left\{ t\in \mathbf{S}(\mathbb{\beth }%
_{j+1}^{-})\ |\ t>\xi \right\} .\ 
\end{equation*}%
Notice that this definition is similar to (\ref{IVcorda}), but in this case $%
x\in \mathbf{No}$ and the relations $<$ and $>$ are in $\mathbf{No}.$ In any
case, $(\mathcal{L}(x),\mathcal{R}(x))$ is a section in $\mathbf{S}(\mathbb{%
\beth }_{j+1}^{-})$ and hence, by Lemma \ref{aria2}, there is a unique $\xi $
between $\mathcal{L}(x)$ and $\mathcal{R}(x)$ and we have that 
\begin{equation*}
i(\xi )=\mathcal{L}(\xi )|\mathcal{R}(\xi )=\mathcal{L}(x)|\mathcal{R}(x)=x
\end{equation*}

$\square $

In conclusion, a section in $\mathbf{S}(\mathbb{\beth }_{j+1}^{-})$
individuates a number $x\in \mathbf{No}$ and a namber in $\xi \in \mathbb{E}$
and these points can be identified. It remains to show that the operations $%
+ $ and $\cdot $ in $\mathbf{No}$ coincide with the operations in $\mathbb{E}
$. We recall that the operations in $\mathbf{No}$ satisfy the following
equalities: 
\begin{eqnarray*}
L(x+y) &=&\left\{ x^{L}+y,\ x+y^{L}\ |\ x^{L}\in L(x),y^{L}\in L(y)\right\}
\\
R(x+y) &=&\left\{ x^{R}+y,\ x+y^{R}\ |\ x^{R}\in L(x),y^{R}\in L(y)\right\}
\end{eqnarray*}%
and%
\begin{eqnarray*}
L(xy) &=&\left\{ x^{L}y+xy^{L}-x^{L}y^{L},x^{R}y+xy^{R}-x^{R}y^{R}\ |\
x^{L}\in L(x),..,y^{R}\in L(y)\right\} \\
R(xy) &=&\left\{ x^{L}y+xy^{R}-x^{L}y^{R},x^{R}y+xy^{L}-x^{R}y^{L}\ |\
x^{L}\in L(x),..,y^{R}\in L(y)\right\}
\end{eqnarray*}

\begin{lemma}
\label{LARA}If $\xi ,\zeta \in \mathbb{R}(\mathbb{\beth }_{j}),$ then 
\begin{equation*}
i(\xi +\zeta )=i(\xi )+i(\zeta )
\end{equation*}%
\begin{equation*}
i(\xi \zeta )=i(\xi )\cdot i(\zeta )
\end{equation*}
\end{lemma}

\textbf{Proof}: Let $x=i(\xi )$ and $y=i(\zeta )$. Since%
\begin{equation*}
\mathcal{L}(\xi )|\mathcal{R}(\xi )=x=L(x)|R(x)
\end{equation*}%
we have that%
\begin{equation*}
I\left( \mathcal{L}(\xi ),\mathcal{R}(\xi )\right) =I(L(x),R(x))
\end{equation*}%
and similarly%
\begin{equation*}
\left( \mathcal{L}(\zeta ),\mathcal{R}(\zeta )\right) =(L(y),R(y))
\end{equation*}%
Then, it is easy to check that%
\begin{equation*}
I\left( \mathcal{L}(\xi +\zeta ),\mathcal{R}(\xi +\zeta )\right)
=I(L(x+y),R(x+y)).
\end{equation*}%
In conclusion%
\begin{equation*}
i(\xi )+i(\zeta )=x+y=L(x+y)|R(x+y)=\mathcal{L}(\xi +\zeta )|\mathcal{R}(\xi
+\zeta )=i(\xi +\zeta )
\end{equation*}%
Using a similar arguments we can prove that $i(\xi )\cdot i(\zeta )=i(\xi
\zeta )$ provided that $\xi ,\zeta >0;$ if $\xi \ $or$\ \zeta $ is not
positive the conclusion follows from standard algebraic manipulations.

$\square $

\bigskip

In conclusion, we have proved the following theorem:

\begin{theorem}
\label{LUNA}The map (\ref{lara}) is a field isomorphism.
\end{theorem}

This result can be resumed as follows:

\begin{corollary}
Let $\mathbf{Ord\ }$be the pseudoclass of ordinal numbers equipped with the
natural operations + and $\cdot \ $and let $\mathbb{F}$ be the smallest real
closed field containing $\mathbf{Ord;}$ then $\mathbb{F}$ is isomorphic to
both $\mathbb{E}$ and $\mathbf{No.}$
\end{corollary}

From now on the sets $\mathbb{E}$ and $\mathbf{No}$ will be identified and
and every Euclidean number will be considered also a surreal number.

\begin{remark}
Th. \ref{LUNA} is new and suggests many directions for developing the theory
of both hyperreal and surreal numbers. For example, given a real function,
we can study the relationship between its natural extension in the set of
hyperreal numbers and its extension in the set of surreal numbers (when and
where it exists). Or, the relationship between hyperfinite and surreal sums.
\end{remark}

\section{A construction of the numerosities\label{CN}}

A numerosity theory is based on a counting system $(\Lambda ,\mathbf{Num,}%
\mathfrak{num})$ which satisfies the Euclid's Principle. Even if the
Euclid's principle is a natural request, it is necessary to prove that it is
consistent with the definition of counting system. We will prove such a
consistency presenting a model based on a peculiar type of labellings $%
\mathfrak{L}$ called \textit{label-tree}.

\subsection{The label-trees\label{LLL}}

Our construction of a numerosity theory is based on a special labelling. In
turn, this labeling is based on a peculiar partial order relation.

\begin{definition}
\label{LT}If $X\subseteq \Lambda ,\ $the triple $(X\mathfrak{,}\sqsubseteq
,\oplus \mathfrak{)}$ is called \textquotedblleft \textbf{pivotal tree}%
\textquotedblright\ if:

\begin{enumerate}
\item $\varnothing \in X$ and $\forall x\in X,\ \varnothing \sqsubseteq x.$

\item $\oplus :X\backslash \{\varnothing \}\rightarrow X\backslash
\{\varnothing \}$ is an injective map; the image $a^{\oplus }$ of an element
will be called successor of \textquotedblleft $a$\textquotedblright ; also,
we will use the notation $x^{m\oplus }:=\left[ x^{(m-1)\oplus }\right]
^{\oplus };$

\item $(X\mathfrak{,}\sqsubseteq \mathfrak{)}$ is a directed set\footnote{$(X%
\mathfrak{,}\sqsubseteq \mathfrak{)}$ is called directed set if $\sqsubseteq 
$ is a preorder relation and 
\begin{equation*}
\forall x,y\in X,\ \exists z\in X,\ \left( x\sqsubseteq z\ \ and\ \
y\sqsubseteq z\right) ;
\end{equation*}%
$\ $As usual, we will employ also the following notation: 
\begin{equation*}
b\sqsupseteq a:\Leftrightarrow a\sqsubseteq b;\ \ \ \ a\equiv
b:\Leftrightarrow \left( a\sqsubseteq b\ \text{and}\ b\sqsubseteq a\right)
.\ 
\end{equation*}%
} such that

\begin{enumerate}
\item $\forall b\in \mathbf{Fin},\ \left( a\in b\ \text{or}\ a\subseteq
b\right) \Rightarrow a\sqsubseteq b;$

\item \label{aaa}$a\sqsubseteq b\Rightarrow \exists m\geq 1,\ a^{m\oplus
}\equiv b;$

\item \label{ccc}$\forall a\in X,$ the set $\{x\in X\ |\ x\sqsubseteq a\}$
is finite.
\end{enumerate}
\end{enumerate}
\end{definition}

\FRAME{ftbpFU}{4.324in}{1.9061in}{0pt}{\Qcb{Schematic representation of a
pivotal tree. Horizontal arrows connect points such that $a\equiv b$.}}{}{%
reticolo7.jpg}{\special{language "Scientific Word";type
"GRAPHIC";maintain-aspect-ratio TRUE;display "USEDEF";valid_file "F";width
4.324in;height 1.9061in;depth 0pt;original-width 20.3952in;original-height
8.9532in;cropleft "0";croptop "1";cropright "1";cropbottom "0";filename
'reticolo7.jpg';file-properties "XNPEU";}}

Let us analyze the structure of a pivotal tree. A pivotal pivotal tree can
be regarded as an oriented graph in which the vertices are the element of $X$
and the oriented arcs have the form $(a,a^{\oplus }).$ Notice that every
chain (by virtue of (\ref{aaa}) and (\ref{ccc})) has a minimum point but not
a maximum. In general $\sqsubseteq $ is a preorder relation and not a
partial order relation, then $(X,\sqsubseteq \mathfrak{)}$ is not a lattice;
however, if we combine $\sqsubseteq $ and $\oplus $, by Def. \ref{LT}-(\ref%
{aaa}), we get a partial order relation%
\begin{equation*}
b\sqsubseteq _{\oplus }a:\Leftrightarrow \exists m\geq 1,\ b^{m\oplus }=a.
\end{equation*}%
\ Then, we can define the \textit{join} $\vee $ and the \textit{meet} with
respect to $\sqsubseteq _{\oplus },$ namely, 
\begin{eqnarray*}
x\vee y &=&\min \{z\ |\ x\sqsubseteq _{\oplus }z\ \text{and}\ \ y\sqsubseteq
_{\oplus }z\}. \\
x\wedge y &=&\max \{z\ |\ z\sqsubseteq _{\oplus }x\ \text{and}\ z\sqsubseteq
_{\oplus }y\}
\end{eqnarray*}%
Notice that, if $x\wedge y\neq \varnothing ,$ then $\exists m\geq 1,\
x^{m\oplus }=y$ or $y^{m\oplus }=x.$

\bigskip

\textbf{Examples}: A trivial example of pivotal tree is given by $(\mathbf{%
Ord}\mathfrak{,}\sqsubseteq ,\oplus \mathfrak{);}$ if we set $\gamma
^{\oplus }=\gamma +1$ and%
\begin{equation*}
b\sqsubseteq a:\Leftrightarrow \exists m\in \mathbb{N},\ b^{m\oplus }=b+m=a.
\end{equation*}

Let see an other example. We set 
\begin{equation}
a\Subset b:\Leftrightarrow \exists b_{1},...,b_{n},\in \mathbf{Fin},\left( \
a\in b_{1}\in b_{2}...\in b_{n}=b\ \ \text{or}\ \ a\subseteq b\right) ,
\label{grulla}
\end{equation}%
moreover, we take a well ordering $\{a_{j}\}_{j\in \mathbf{Ord}}\ $of $%
\Lambda $ consistent with $\Subset ,$ namely%
\begin{equation*}
a_{j}\Subset a_{k}\Rightarrow j<k
\end{equation*}%
We define $\oplus $ and $\sqsubseteq $ as follows, 
\begin{equation}
a_{j}\sqsubset a_{k}:\Leftrightarrow \left( j<k\ \ \text{and}\ \
a_{j}\Subset a_{k}\right) ;\ \ \ a_{j}\equiv a_{k}:\Leftrightarrow
a_{j}=a_{k}  \label{exa}
\end{equation}%
\begin{equation*}
a_{j}^{\oplus }:=a_{k}\ \ \text{where\ \ }\ k=\min \{p\in \mathbf{Ord\ }|\
a_{p}\sqsupseteq a_{k}\}
\end{equation*}%
then it is easy to check that $(\Lambda \mathfrak{,}\sqsubseteq ,\ \oplus 
\mathfrak{)}$ is a pivotal tree.

\begin{theorem}
\label{PP}Given a pivotal tree $(\Lambda \mathfrak{,\sqsubseteq ,\oplus ),}$
there is a label lattice $(\mathfrak{L,\subseteq )}$, called \textbf{%
label-tree,} and a label map%
\begin{equation*}
\ell :\Lambda \rightarrow \mathfrak{L,\ \ }\ell (a):=\dbigcup \{x\in \Lambda
\ |\ x\sqsubseteq a\}=\{x\in \Lambda \ |\ \exists b\subseteq a,\
x\sqsubseteq b\}
\end{equation*}%
such that

\begin{enumerate}
\item \label{f1}if $\lambda ,\mu \in \mathfrak{L}$,\ then $\lambda \cap \mu
=\lambda \wedge \mu \in \mathfrak{L\ }$and $\lambda \vee \mu \in \mathfrak{L}%
;$

\item \label{q2}$a\sqsubseteq b\Rightarrow \ell (a)\subseteq \ell (b);$

\item \label{q}$\ell \left( \ell (a)\right) =\ell (a);$

\item \label{f2}$\ell (a\wedge b)=\ell (a)\wedge \ell (b),\ $

\item \label{f3}$\ell (a\vee b)=\ell (a)\vee \ell (b),\ $

\item \label{q4}$\ell (\left\{ a,b\right\} )=\ell (\{a\})\vee \ell (\{b\});$

\item \label{q1}$\ell (\left( a,b\right) )=\ell (\{\{a\}\})\vee \ell
(\{\{b\}\}).$
\end{enumerate}
\end{theorem}

\textbf{Proof}: Given $a\in \Lambda \mathfrak{,}$ we set%
\begin{equation*}
\lambda _{a}:=\dbigcup \{x\in \Lambda \ |\ x\sqsubseteq a\}
\end{equation*}%
and $\mathfrak{L}=\{\lambda _{a}\ \ |\ a\in \Lambda \}.$

First of all we have to prove that $\mathfrak{L}$ is a labelling, namely
that $(\mathfrak{L,\subseteq )}$ the requests of Def. \ref{LU} are
satisfied. (\ref{a}) and (\ref{c}) are immediate. Let us see (\ref{b}). We
have that 
\begin{eqnarray}
\lambda _{a}\cap \lambda _{b} &=&\left( \dbigcup \{x\in \Lambda \ |\
x\sqsubseteq a\}\right) \cap \left( \dbigcup \{x\in \Lambda \ |\
x\sqsubseteq b\}\right)  \notag \\
&=&\dbigcup \{x\in \Lambda \ |\ x\sqsubseteq a\cap b\}=\lambda _{a\cap b}\in 
\mathfrak{L}  \label{f10}
\end{eqnarray}%
Similarly, 
\begin{eqnarray}
\lambda _{a}\cup \lambda _{b} &=&\left( \dbigcup \{x\in \Lambda \ |\
x\sqsubseteq a\}\right) \cup \left( \dbigcup \{x\in \Lambda \ |\
x\sqsubseteq b\}\right)  \notag \\
&=&\dbigcup \{x\in \Lambda \ |\ x\sqsubseteq a\ \text{or }x\sqsubseteq \cup
b\}  \label{f11} \\
&\subseteq &\dbigcup \{x\in \Lambda \ |\ x\sqsubseteq _{\oplus }a\vee
b\}=\lambda _{a\mathfrak{\vee }b}\mathfrak{.}  \notag
\end{eqnarray}

Now, let us prove that 
\begin{equation*}
\ell (a):=\dbigcap \left\{ \mu \in \mathfrak{L}\ |\ a\in \mu \right\}
=\lambda _{a}
\end{equation*}%
It holds%
\begin{equation*}
\ell (a)=\dbigcap \left\{ \lambda _{x}\ |\ x\in \Lambda ,\ a\in \lambda
_{x}\right\}
\end{equation*}%
Since $a\in \lambda _{a},$ we have that $\ell (a)\subseteq \lambda _{a}$. If 
$a\in \mu ,$ and $b\sqsubseteq a,$ then $b\sqsubseteq \mu ;$ therefore $b\in
\lambda _{a}\Rightarrow b\in \mu .$ Then%
\begin{equation*}
\lambda _{a}\subseteq \dbigcap \left\{ \mu \in \mathfrak{L}\ |\ a\in \mu
\right\} =\ell (a).
\end{equation*}

(\ref{f1}) follows from (\ref{f10}) and (\ref{f11})

(\ref{q2}) - Trivial.

(\ref{q}) - We have that 
\begin{equation*}
\ell (a)=\{x\in \Lambda \ |\ x\sqsubseteq a\}=\{x\in \Lambda \ |\
x\sqsubseteq \lambda _{a}\}=\ell (\lambda _{a})=\ell \left( \ell (a)\right) .
\end{equation*}

(\ref{f2}) - We have that 
\begin{eqnarray*}
\ell (a)\wedge \ell (b) &=&\lambda _{a}\cap \lambda _{b}=\dbigcup \{x\in
\Lambda \ |\ x\sqsubseteq a\ \text{and\ }x\sqsubseteq b\} \\
&=&\dbigcup \{x\in \Lambda \ |\ x\sqsubseteq a\wedge b\}=\ell (a\wedge b)\ 
\end{eqnarray*}

(\ref{f3}) - We have that 
\begin{eqnarray*}
\ell (a)\vee \ell (b) &=&\lambda _{a}\vee \lambda _{b}=\{x\in \Lambda \ |\
x\sqsubseteq a\ \text{or\ }x\sqsubseteq b\} \\
&=&\dbigcup \{x\in \Lambda \ |\ x\sqsubseteq a\vee b\}=\ell (a\vee b)\ 
\end{eqnarray*}

(\ref{q4})- By Def.\ref{LT}-(\ref{a}), we have that%
\begin{equation*}
\left\{ a\right\} \sqsubseteq \left\{ a,b\right\} ,\ \left\{ b\right\}
\sqsubseteq \left\{ a,b\right\}
\end{equation*}%
and hence 
\begin{equation*}
\left\{ a\right\} \vee \left\{ b\right\} \sqsubseteq \left\{ a,b\right\}
=\left\{ a\right\} \cup \left\{ b\right\} \sqsubseteq \left\{ a\right\} \vee
\left\{ b\right\} ;
\end{equation*}%
then $\left\{ a\right\} \vee \left\{ b\right\} =\left\{ a,b\right\} \ $and%
\begin{equation*}
\ell (\left\{ a,b\right\} )=\ell (\left\{ a\}\vee \{b\right\} )=\ell
(\{a\})\vee \ell (\{b\})
\end{equation*}

(\ref{q1}) - By (\ref{q4}), identifying $\left( a,b\right) $ with a
Kuratosky pair, 
\begin{eqnarray*}
\ell (\left( a,b\right) ) &=&\ell (\left\{ \left\{ a\right\} ,\left\{
a,b\right\} \right\} )=\ell (\{\left\{ a\right\} \})\vee \ell (\{\left\{
a,b\right\} \}) \\
&=&\ell (\{\left\{ a\right\} \})\vee \ell (\{\left\{ a\right\} \})\vee \ell
(\{\left\{ b\right\} \})=\ell (\{\left\{ a\right\} \})\vee \ell (\{\left\{
b\right\} \}).
\end{eqnarray*}

$\square $

\bigskip

By this theorem, if we restrict $(\Lambda \mathfrak{,\mathfrak{\vee ,\wedge }%
)}$ to $\mathfrak{L,}$ we get a sublattice $(\mathfrak{L}\mathfrak{,%
\mathfrak{\vee ,\wedge })}$. The label map%
\begin{equation*}
\ell :\Lambda \rightarrow \mathfrak{L}
\end{equation*}%
is an homomorphism from the pivotal tree $(\Lambda \mathfrak{,\mathfrak{\vee
,\wedge })}$ to the label-tree $(\mathfrak{L}\mathfrak{,\mathfrak{\vee
,\wedge }).}$

From the point of view of the graph theory, a label $\ell (a)$ is the union
of all the chains which start at $\varnothing $ and end to $a.$

\textbf{Examples}: In order to familiarize the reader with the notion of
label tree, we will give some trivial example using the partial order
relation (\ref{exa}) with the further assumption that $\forall n\in \mathbb{N%
},\ a_{n}=n$

\begin{itemize}
\item $\ell (3)=\{3,\varnothing \};$

\item $\ell (\mathbb{N})=\{\mathbb{N},\varnothing \};$

\item $\ell (\{2\})=\{\{2\},2,\varnothing \};$

\item $\ell (\{1,2\})=\{\{1,2\},\{1\},\{2\},2,1,\varnothing \}.$
\end{itemize}

\subsection{Basic properties of label-trees}

In this section we will prove some technical lemmas which will be used in
the next sections.

\begin{lemma}
\label{v33}If $\left( \mathfrak{L},\subseteq \right) $ is a label-tree and $%
\lambda ,\mu \in \mathfrak{L}$, then%
\begin{equation*}
\lambda \mathfrak{\mathfrak{\wedge }}\mu =\lambda \ \text{or\ }\lambda 
\mathfrak{\mathfrak{\wedge }}\mu =\mu \ \ \text{or \ }\lambda \mathfrak{%
\mathfrak{\wedge }}\mu =\varnothing .
\end{equation*}
\end{lemma}

\textbf{Proof}: Set $\lambda =\ell (a)$ and $\mu =\ell (b)$. If $%
a\sqsubseteq b$,\ then $\ell (a)\subseteq \ell (b);\ $hence%
\begin{equation*}
\lambda \mathfrak{\mathfrak{\wedge }}\mu =\lambda \mathfrak{\mathfrak{\cap }}%
\mu =\ell (a)\cap \ell (b)=\ell (b)=\mu ;
\end{equation*}%
if $b\sqsubseteq a$,\ we argue in the same way. If $a\ $and $b$ are not
comparable, we claim that $\ell (a)\cap \ell (b)=\{\varnothing \}$; in fact,
if it would exist $c\in \ell (a)\cap \ell (b),\ $we$\ $would have that $%
a=c^{m\oplus }$\ and\ $b=c^{n\oplus };$ namely $a$\ and\ $b$ would be
comparable.

$\square $

\begin{lemma}
\label{puppa}If $\left( \mathfrak{L},\subseteq \right) $ is a label-tree,
there is a well ordering $\left\{ \lambda _{j}\right\} _{j\in \Delta }\ $of $%
\mathfrak{L}\ $such that%
\begin{equation*}
\lambda _{k}\subset \lambda _{j}\Rightarrow k<j.
\end{equation*}
\end{lemma}

\textbf{Proof: }We set $F_{0}=\mathbf{Ato\cup Inf}$ and, for $n\geq 0,$ we
will denote by $F_{n+1}$ the family of elements in $\Lambda $ which have a
"minimal predecessor" in $F_{n},$ namely%
\begin{equation*}
a\in F_{n+1}:\Leftrightarrow \left[ \left( \exists x\in F_{n},\ x^{\oplus
}=a\right) \ \text{and}\ \left( y^{\oplus }=a\ \Rightarrow y\in F_{m},\
m>n\right) \right] .
\end{equation*}%
Next, we take a well-ordering $\{a_{j}\}_{j\in \triangle }$ of $\Lambda
\backslash \{\varnothing \}$ consistent with the sets $F_{n}$ namely if%
\begin{equation*}
\left[ a_{k}\in F_{n}\ \text{and\ }a_{j}\in F_{m},\ m>n\right] \Rightarrow
k<j.
\end{equation*}

Then, $\left\{ \lambda _{j}\right\} _{j\in \Delta }$ satisfies our request.
In fact if $\lambda _{k}=\ell (a_{k})\subset \lambda _{j}=\ell (a_{j}),$ we
have that $a_{k}\in \ell (a_{j})$ and hence $a_{j}\equiv a_{k}^{m\oplus }.$
Then $k<j.$

$\square $

\begin{lemma}
\label{v3}Let $\left\{ \lambda _{j}\right\} _{j\in \Delta }\ $be of a well
ordering $\mathfrak{L}\ $as in lemma \ref{puppa}$.\ $Then, there exists
finite sequence of sets $\{S_{j_{0}},S_{j_{1}}..,S_{j_{n}}\}$ such that%
\begin{equation}
p\neq q\Rightarrow S_{j_{p}}\cap S_{j_{q}}=\varnothing  \label{v4}
\end{equation}%
and 
\begin{equation*}
\lambda _{j}=S_{j_{0}}\cup ...\cup S_{j_{n}}
\end{equation*}
\end{lemma}

\textbf{Proof}: We set $j_{0}=j$ and 
\begin{equation*}
j_{1}:=\max \{k\ |\ \lambda _{k}\subset \lambda _{j_{0}}\},\ \
S_{j_{0}}=\lambda _{j_{0}}\backslash \lambda _{j_{1}};
\end{equation*}%
then, if $\lambda _{j_{1}}\neq \varnothing ,$ 
\begin{equation*}
\lambda _{j}=S_{j_{0}}\cup \lambda _{j_{1}}
\end{equation*}%
If we iterate this operation with $\lambda _{j_{1}},$ we get 
\begin{equation*}
\lambda _{j}=S_{j_{0}}\cup S_{j_{1}}\cup \lambda _{j_{2}}
\end{equation*}%
This process ends when $\{k\ |\ \lambda _{k}\subset \lambda
_{j_{n}}\}=\varnothing .$

$\square $

\subsection{A numerosity counting system}

In this section we will construct a \textit{numerosity counting system}
exploiting a \textit{fine ultrafilter} over the label-tree.

\begin{definition}
Given a label-tree, $(\mathfrak{L,\subseteq )}$ a \textbf{fine ultrafilter\ }
$\mathcal{U}$ over $\mathfrak{L}$ is a family of subsets of $\mathfrak{L}$
which satisfies the following properties:

\begin{enumerate}
\item $\mathfrak{L}\in \mathcal{U}$

\item if $Q\in \mathcal{U}$ and $P\supset Q$,\ then $P\in \mathcal{U}$,

\item if $P,Q\in \mathcal{U}$,\ then $P\cap Q\in \mathcal{U}$,

\item if $Q\in \mathcal{U}$,\ then $\mathfrak{\wp }_{\omega }(\Lambda
)\backslash Q\notin \mathcal{U},$

\item $\forall \lambda \in \mathfrak{L}$,\ $C\left[ \lambda \right] \in 
\mathcal{U};$ here $C\left[ \lambda \right] $ denotes the cone with the
vertex in $\lambda ,$ namely 
\begin{equation}
C\left[ \lambda \right] :=\left\{ \mu \in \mathfrak{L\ }|\ \lambda \subseteq
\mu \right\} .
\end{equation}
\end{enumerate}
\end{definition}

It is well known that the existence of such ultrafilter is a consequence of
Zorn's lemma. As usual, a set $Q\in \mathcal{U}$ is called \textbf{qualified}%
.

\begin{definition}
\label{fighissima}We set%
\begin{equation*}
A\preceq _{\mathfrak{n}}B
\end{equation*}%
if there exists a qualified set $Q\in \mathcal{U}$ such that $\forall
\lambda \in Q$%
\begin{equation}
\left\vert A\cap \lambda \right\vert \leq \left\vert B\cap \lambda
\right\vert .  \label{billo}
\end{equation}
\end{definition}

The following result holds.

\begin{theorem}
\label{T1+}The couple $(\Lambda ,\preceq _{\mathfrak{n}})$ is a comparison
system
\end{theorem}

\textbf{Proof}: First of all, let us see that $\preceq _{\mathfrak{n}}$ is a
preorder relation; if $A\preceq _{\mathfrak{n}}B$ and $B\preceq _{\mathfrak{n%
}}C$, there are two qualified sets $Q_{1},Q_{2}\in \mathcal{U}$ such that%
\begin{equation*}
\forall \lambda \in Q_{1},\ \left\vert A\cap \lambda \right\vert \leq
\left\vert B\cap \lambda \right\vert \ \ \ and\ \ \ \forall \lambda \in
Q_{2},\ \left\vert B\cap \lambda \right\vert \leq \left\vert C\cap \lambda
\right\vert
\end{equation*}%
then, 
\begin{equation*}
\forall \lambda \in Q_{1}\cap Q_{2},\ \left\vert A\cap \lambda \right\vert
\leq \left\vert C\cap \lambda \right\vert ;
\end{equation*}%
since $Q_{1}\cap Q_{2}\in \mathcal{U}$, $A\preceq _{\mathfrak{n}}C.$

Now let us prove the points (\ref{C1}) - (\ref{C2}) of Def. \ref{C}.

\ref{C}-(\ref{C1}) - Trivial.

\ref{C}-(\ref{C3}) - \ If $A\cap B=A^{\prime }\cap B^{\prime }=\varnothing
,\ $and $A\cong A^{\prime },\ B\cong B^{\prime },$ then there are two
qualified sets $Q_{1},Q_{2}\in \mathcal{U}$ such that%
\begin{equation*}
\forall \lambda \in Q_{1},\ \left\vert A\cap \lambda \right\vert =\left\vert
A^{\prime }\cap \lambda \right\vert ,\ \forall \lambda \in Q_{2},\
\left\vert B\cap \lambda \right\vert =\left\vert B^{\prime }\cap \lambda
\right\vert ,
\end{equation*}%
then, $\forall \lambda \in Q_{1}\cap Q_{2},\ $ 
\begin{eqnarray*}
\left\vert \left( A\cup B\right) \cap \lambda \right\vert &=&\ \left\vert
\left( A\cap \lambda \right) \cup \left( B\cap \lambda \right) \right\vert
=\left\vert A\cap \lambda \right\vert +\left\vert B\cap \lambda \right\vert
\\
&=&\left\vert A^{\prime }\cap \lambda \right\vert +\left\vert B^{\prime
}\cap \lambda \right\vert =\ \left\vert \left( A^{\prime }\cup B^{\prime
}\right) \cap \lambda \right\vert .
\end{eqnarray*}

\ref{C}-(\ref{C4}) - Let 
\begin{equation}
C\left[ \ell ((a,b))\right] :=\left\{ \lambda \in \mathfrak{L\ }|\ \lambda
\supseteq \ell ((a,b))\right\}  \label{cono}
\end{equation}%
denote the \textit{cone} over $\ell ((a,b))$. By Prop.\ref{PP} (\ref{q1}%
)-(7), and the request of Def. \ref{C}-(\ref{C4}), $\forall \lambda \in C%
\left[ \ell ((a,b))\right] ,$ we have that, 
\begin{eqnarray*}
\left( A\times B\right) \cap \lambda &=&\left\{ x\in A\times B\mathfrak{\ }%
|\ x\in \lambda \right\} \\
&=&\left\{ \left( a,b\right) \mathfrak{\ }|\ a\in A,\ b\in B\ \text{and}\
\left( a,b\right) \in \lambda \right\} \\
&=&\left\{ \left( a,b\right) \mathfrak{\ }|\ a\in A,\ b\in B\ \text{and}\
\ell (\left( a,b\right) )\in \lambda \ \right\} \\
&=&\left\{ \left( a,b\right) \mathfrak{\ }|\ a\in A,b\in B\ \text{and}\ \ell
(\{\{a\}\})\vee \ell (\{\{b\}\})\in \lambda \ \right\} \\
&=&\left\{ a\in A\mathfrak{\ }|\ \ell (\{\{a\}\})\in \lambda \ \right\}
\times \left\{ b\in B\mathfrak{\ }|\ \ell (\{\{b\}\})\in \lambda \ \right\}
\end{eqnarray*}%
\ Since $\ell (a)\subseteq \ell (\{\{a\}\}),\ $we have that $\left( a\in A,\
\ell (\{\{a\}\})=\lambda \Leftrightarrow \ell (a)\subseteq \lambda \right) ,$
then, 
\begin{eqnarray*}
\left( A\times B\right) \cap \lambda &=&\left\{ a\in A\mathfrak{\ }|\ \ell
(a)\in \lambda \right\} \times \left\{ b\in B\mathfrak{\ }|\ \ell (b)\in
\lambda \right\} \\
&=&\left( A\cap \lambda \right) \times \left( B\cap \lambda \right) .
\end{eqnarray*}%
In conclusion, $\forall \lambda \in \mathfrak{L},$ 
\begin{equation*}
\left\vert \left( A\times B\right) \cap \lambda \right\vert =\left\vert
\left( A\cap \lambda \right) \times \left( B\cap \lambda \right) \right\vert
=\left\vert A\cap \lambda \right\vert \cdot \left\vert B\cap \lambda
\right\vert
\end{equation*}

\ref{C}-(\ref{C5}) - Let $C\left[ \mathbb{\ell }(a)\right] $ denotes the
cone defined as in \ref{cono}. By (\ref{C4}), for every label $\lambda \in C%
\left[ \mathbb{\ell }(a)\right] $ 
\begin{equation*}
\left\vert \left( \{a\}\times B\right) \cap \lambda \right\vert =\left\vert
\{a\}\cap \lambda \right\vert \cdot \left\vert B\cap \lambda \right\vert
=1\cdot \left\vert B\cap \lambda \right\vert =\left\vert B\cap \lambda
\right\vert .
\end{equation*}%
The conclusion follows by the fact that $C\left[ \mathbb{\ell }(a)\right] $
is qualified.

\ref{C}-(\ref{C2}) - If $A\preceq _{\mathfrak{n}}B,$ $\exists Q\in \mathcal{U%
}$ such that $\forall \lambda \in Q$%
\begin{equation}
\left\vert A\cap \lambda \right\vert \leq \left\vert B\cap \lambda
\right\vert .  \label{pli}
\end{equation}%
Let $\left\{ \lambda _{j}\right\} _{j\in \Delta }\ $be of a well ordering $%
Q\ $as in lemma \ref{puppa}. By lemma \ref{v3}, we have that $\lambda
_{j}=S_{j_{0}}\cup ...\cup S_{j_{n}}$ with $j_{0}=j.$ We set $%
X_{j_{m}}:=A\cap S_{j_{m}}$ ($0<m<n$) and hence%
\begin{equation*}
A\cap \lambda _{j}=X_{j_{0}}\cup X_{j_{1}}\cup ...\cup X_{j_{n}}
\end{equation*}%
with $j_{0}>...>j_{n}$ and by (\ref{v4}) 
\begin{equation*}
j_{h}\neq j_{m}\ \Leftrightarrow X_{j_{h}}\cap X_{j_{m}}=\varnothing .
\end{equation*}%
Put $j_{0}=\min \left( \Delta \right) .$ Since $\left\vert A\cap \lambda
_{j_{0}}\right\vert \leq \left\vert B\cap \lambda _{j_{0}}\right\vert ,$ $%
\left\vert X_{j_{0}}\right\vert <\left\vert B\cap \lambda
_{j_{0}}\right\vert ;$ then we can take $Y_{j_{0}}\subset B\cap \lambda
_{j_{0}}$ such that $\left\vert Y_{j_{0}}\right\vert =\left\vert
X_{j_{0}}\right\vert $. Now, we claim that for every $k\leq j,$ there are
sets $Y_{k}\subset B$ such that%
\begin{equation}
Y_{k}\subseteq \left( B\cap \lambda _{k}\right) \backslash \left[
Y_{k_{0}}\cup ...\cup Y_{k_{0}}\right] ;\ \ \ \left\vert Y_{k}\right\vert =\
\left\vert X_{k}\right\vert ;\   \label{plu}
\end{equation}%
and%
\begin{equation*}
j_{n}\neq j_{m}\ \Leftrightarrow Y_{j_{n}}\cap Y_{j_{m}}=\varnothing .
\end{equation*}

We argue by induction over $j\in \Delta .$ If $j_{0}=\min \left( \Delta
\right) ,$ (\ref{plu}) holds by the definition of $Y_{j_{0}}.$ Now, we
assume that (\ref{plu}) holds for every $k<j$, then 
\begin{eqnarray*}
\left\vert X_{j}\right\vert &=&\left\vert \left( A\cap \lambda _{j}\right)
\backslash \left[ X_{j_{n}}\cup ...\cup X_{j_{0}}\right] \right\vert
=\left\vert A\cap \lambda _{j}\right\vert -\left\vert X_{j_{0}}\right\vert
-\left\vert X_{j_{1}}\right\vert -...-\left\vert X_{j_{n}}\right\vert \\
&\leq &\left\vert B\cap \lambda _{j}\right\vert -\left\vert
Y_{j_{0}}\right\vert -\left\vert Y_{j_{1}}\right\vert -...-\left\vert
Y_{j_{n}}\right\vert =\left\vert \left( B\cap \lambda _{j}\right) \backslash 
\left[ Y_{j_{0}}\cup ...\cup Y_{j_{n}}\right] \right\vert
\end{eqnarray*}%
Hence, it is possible to take a set $Y_{j}\subseteq \left( B\cap \lambda
_{j}\right) \backslash \left[ Y_{j_{1}}\cup ...\cup Y_{j_{n}}\right] $ such
that $\left\vert Y_{j}\right\vert =\left\vert X_{j}\right\vert .$ Hence (\ref%
{plu}) holds for every $j\in \triangle $. Finally, we set%
\begin{equation*}
A^{\prime }:=\dbigcup\limits_{j\in \triangle }Y_{j}
\end{equation*}%
and we have that 
\begin{equation*}
\left\vert A_{j}^{\prime }\cap \lambda _{j}\right\vert =\left\vert Y_{j}\cup
Y_{j_{1}}\cup ...\cup Y_{j_{n}}\right\vert =\left\vert X_{j}\cup
X_{j_{1}}\cup ...\cup X_{j_{n}}\right\vert =\left\vert A_{j}\cap \lambda
_{j}\right\vert
\end{equation*}

$\square $

\begin{corollary}
If we set%
\begin{equation*}
\mathfrak{num}(A):=\Phi \left( \left[ A\right] _{\cong _{\mathfrak{n}%
}}\right) .
\end{equation*}%
and 
\begin{equation*}
\mathbf{Num:=}\left\{ x\in \mathbb{E}\ |\ \exists A\in \mathcal{\Lambda },\
x=\mathfrak{num}(A)\right\}
\end{equation*}%
then $(\Lambda _{\flat },\mathbf{Num,}\mathfrak{num})$ is a numerosity
theory.
\end{corollary}

\textbf{Proof}: It is immediate to see that $(\Lambda ,\mathbf{Num,}%
\mathfrak{num})$ satisfies the Euclid's Principle.

$\square $

\section{Special properties of numerosities\label{SP}}

The properties of numerosities described in the previous sections are shared
by every numerosity counting system. However if we want to answer to some
specific questions, the information given by Def. \ref{C} is not sufficient.
For example, we cannot compare the numerosity of $\wp (\mathbb{N})$ and the
numerosity of $\mathbb{R}$. In general given two sets $A$ and $B$, not
always it is possible to compare the $\Lambda $-limits of $|A\cap \lambda |$
and $|B\cap \lambda |$ since we do not have enough information. Actually
different answers are consistent with the definitions/axioms of the theory.
Therefore we can add new axioms and to check that they are consistent.
However, from a technical point of view, it is easier to take an appropriate
label-tree $\mathfrak{L}$ so that we can compare $|A\cap \lambda |$ and $%
|B\cap \lambda |$ for every sufficiently large label $\lambda \in \mathfrak{L%
}$. In other words, rather then adding new axioms, we can choose a suitable
label-tree. Of course, this choice is absolutely arbitrary since it
correspond to the addiction of independent axioms. Anyway, we can enrich the
theory with new properties which we will call "special".

\subsection{The general strategy\label{LL}}

\bigskip

In order to define appropriate label-trees, we will exploit the following
theorem:

\begin{theorem}
\label{laura}Given a set $\mathfrak{S}\subset \wp _{\omega }(\Lambda ),$ if $%
\{\sigma \in \mathfrak{S}\ |\ a\Subset \sigma \}\neq \varnothing $, we put%
\begin{equation*}
\ell _{\mathfrak{S}}(a)=\dbigcap \{\sigma \in \mathfrak{S}\ |\ a\Subset
\sigma \};
\end{equation*}%
then, there exist a label-tree $(\mathfrak{L}_{\mathfrak{S}}\mathfrak{%
,\subseteq )}$ and a fine ultrafilter over $\mathfrak{L}_{\mathfrak{S}}$,
such that 
\begin{equation}
\ell _{\mathfrak{S}}(a)\in \mathfrak{L.}  \label{78}
\end{equation}%
and the set%
\begin{equation*}
Q_{\mathfrak{S}}:=\{\ell (a)\ |\ a\in \dbigcup \mathfrak{S}\}
\end{equation*}%
is qualified relatively\footnote{%
Given an ultrafilter $\mathcal{U}$ over $\mathfrak{L}$ and set $X,$ we say
that a set $Q$ is qualified relatively to $X$ if there exists a set $\hat{Q}%
\in \mathcal{U}$ such that%
\begin{equation*}
Q=\{\hat{\lambda}\cap X\ |\ \hat{\lambda}\in \hat{Q}\}
\end{equation*}%
} to $\dbigcup \mathfrak{S.\ }$We will refer to $(\mathfrak{L}_{\mathfrak{S}}%
\mathfrak{,\subseteq )}$ as to a label tree \textbf{induced} by $\mathfrak{S}%
\mathfrak{.}$
\end{theorem}

\textbf{Proof}: First of all we define a preorder relation over $\Lambda $
as follows: 
\begin{equation}
a\sqsubseteq \mathfrak{_{\mathfrak{S}}}b:\Leftrightarrow \left( \ell _{%
\mathfrak{S}}(a)\subseteq \ell _{\mathfrak{S}}(b)\ \ \text{or\ \ }a\Subset
b\right)  \label{ble}
\end{equation}%
Next we define a map $\oplus :\Lambda \rightarrow \Lambda $ so that $%
(\Lambda \mathfrak{,}\sqsubseteq _{\mathfrak{S}},\oplus \mathfrak{)}$ be a
pivotal tree as follows: we set%
\begin{eqnarray*}
F_{0} &=&\mathbf{Ato}\cup \mathbf{Inf}\cup \{\varnothing \} \\
F_{n+1} &=&\wp _{n}(F_{n})\cup F_{n}
\end{eqnarray*}%
where $\wp _{n}(X):=\{x\in \wp (X)\ |\ \left\vert x\right\vert \leq n\}$.
Then $\Lambda =\dbigcup\limits_{n\in \mathbb{N}}F_{n}.$ Let $\{a_{j}\}_{j\in 
\mathbf{Ord}}$ be a well ordering of $\Lambda $ consistent with the $F_{n}$%
's, namely such that, $a_{0}=\varnothing $ and 
\begin{equation*}
\left( a_{j}\in F_{n}\ \ \text{and}\ \ a_{h}\in F_{m}\backslash F_{n}\right)
\Rightarrow j<k.
\end{equation*}%
Now we set, for $j\geq 1$, 
\begin{equation*}
a_{j}^{\oplus }=a_{k};\ \ k=\min \{q\ |\ k>j,\ a_{j}\sqsubseteq \mathfrak{_{%
\mathfrak{S}}}a_{k}\}
\end{equation*}%
Notice that this definition is well posed, in fact $a_{j}\sqsubseteq _{%
\mathfrak{S}}\{a_{j}\}$ and since $\{a_{j}\}=a_{p}$ for some $p>j,$ we have
that $\{q\ |\ k>j,\ a_{j}\sqsubseteq _{\mathfrak{S}}a_{k}\}\neq \varnothing $%
.\ It is easy to see that $(\Lambda \mathfrak{,}\sqsubseteq _{\mathfrak{S}%
},\oplus \mathfrak{)}$ is a pivotal tree and hence by Th. \ref{PP} there
exists a label-tree $(\mathfrak{L}_{\mathfrak{S}},\subseteq ).$\ (\ref{78})
is satisfied by construction. Now, it is sufficirnt to take a fine
ultrafilter over $\mathfrak{L}_{\mathfrak{S}}$ which contains $Q_{\mathfrak{S%
}}.$

$\square $

\bigskip

The introduction of $\mathfrak{S}$ and Th. \ref{laura} allows to compare the
numerosities of suitable set exploiting the following proposition:

\begin{proposition}
\label{PR}If $A,B\subset \dbigcup \mathfrak{S,}$ and $\forall \lambda =Q_{%
\mathfrak{S}},\ \left\vert A\cap \lambda \right\vert =\left\vert B\cap
\lambda \right\vert \ $then $\mathfrak{num}(A)=\mathfrak{num}(B)\mathfrak{.}$
\end{proposition}

\textbf{Proof}: Since $Q_{\mathfrak{S}}$ is qualified relatively to $%
\dbigcup \mathfrak{S,}$ there exists a qualified set $\hat{Q}\in \mathcal{U}$%
; thus $\forall \hat{\lambda}\in \hat{Q},$ we have that$\ \left\vert A\cap 
\hat{\lambda}\right\vert =\ \left\vert A\cap \lambda \right\vert $ and $%
\left\vert B\cap \hat{\lambda}\right\vert =\ \left\vert B\cap \lambda
\right\vert .$ Hence $\forall \hat{\lambda}\in \hat{Q},\ \left\vert A\cap
\lambda \right\vert =\left\vert B\cap \lambda \right\vert ;\ $the conclusion
follows by taking the $\Lambda $-limit.

$\square $

\bigskip

Clearly, if we expand $\mathfrak{S,}$ we get more information. Of course, it
is necessary to do it in a suitable way.

\begin{definition}
We say that a set $\mathfrak{S}_{2}$ is \textbf{compatible }with $\mathfrak{S%
}_{1}$, if $\forall a,b\in \mathfrak{S}_{1}\cap \mathfrak{S}_{2},$%
\begin{equation*}
a\cap b\in \mathfrak{S}_{1}
\end{equation*}
\end{definition}

If $\mathfrak{S}_{1}$ and $\mathfrak{S}_{2}$ are compatible, we can expand $%
\mathfrak{S}_{1}$ by putting $\mathfrak{S}=\mathfrak{S}_{1}\cup \mathfrak{S}%
_{2}$. The compatibility guarantees that for every $a\subset \dbigcup 
\mathfrak{S}_{1}$ 
\begin{equation*}
\ell _{\mathfrak{S}}(a)=\ell _{\mathfrak{S}_{1}}(a).
\end{equation*}

\subsection{The numerosities of some subsets of$\ \mathbb{R}$}

\subsubsection{Numerosity of the natural numbers}

Our goal is to define a label-tree that provides "nice" properties to the
subset of $\mathbb{N}$ following the procedure described in section \ref{LL}%
. We set 
\begin{equation*}
\mathfrak{S}\left( \mathbb{N}\right) :=\left\{ \{0,...,n\}\ |\ n\in \Pi
\right\} \}\ \text{where}\ \ \Pi :=\left\{ m!^{m!}\ |\ m\in \mathbb{N}%
\right\} \ \ 
\end{equation*}%
and\ we apply Th. \ref{laura} to induce the label-tree $\left( \mathfrak{L}_{%
\mathbb{N}},\subseteq \right) $. With this choice we have that $\forall n\in 
\mathbb{N}$,%
\begin{equation*}
\ell _{_{\mathbb{N}}}(n)=\ell _{_{\mathbb{N}}}(\{n\})=\{0,1,2,...,m!^{m!}\}
\end{equation*}%
where $m=\min \{m\in \mathbb{N}\ |\ m!^{m!}\geq n\}.$

In order to simplify some algebraic manipulations, it is useful to introduce
the number 
\begin{equation}
\mathbf{\alpha }:=\mathfrak{num}(\mathbb{N}^{+})=\lim_{\lambda \uparrow
\Lambda }\ \left\vert \mathbb{N}^{+}\cap \lambda \right\vert =\omega -1
\label{alfa0}
\end{equation}

The reason for such a choice of $\mathfrak{S}\left( \mathbb{N}\right) $ is
to ensure the following properties of $\mathbf{\alpha }$:

\begin{theorem}
\label{basic}Let $p\in \mathbb{N}^{+}$. Then

\begin{enumerate}
\item if for $i=0,\dots ,n-1$%
\begin{equation*}
A_{i}=\{n\in \mathbb{N}^{+}\mid n\equiv i\ \ \func{mod}p\}.
\end{equation*}%
then $\mathfrak{num}(A_{i})=\frac{\mathbf{\alpha }}{p}$;

\item if%
\begin{equation*}
E_{p}=\{x^{p}\in \mathbb{N}^{+}\mid x\in \mathbb{N}^{+}\}
\end{equation*}%
then%
\begin{equation*}
\mathfrak{num}(E_{p})=\mathbf{\alpha }^{\frac{1}{p}}.
\end{equation*}
\end{enumerate}
\end{theorem}

\textbf{Proof}: (1) - For every $\lambda =\{0,1,\dots ,m!^{m!}\}\in Q_{%
\mathfrak{S}(\mathbb{N})}$, $m\geq p,$ we have that 
\begin{equation*}
|A_{i}\cap \lambda |=\frac{m!^{m!}}{p}=\frac{\left\vert \mathbb{N}^{+}\cap
\lambda \right\vert }{p}
\end{equation*}%
Hence,%
\begin{equation*}
\mathfrak{num}(A_{i})=\lim_{\lambda \uparrow \Lambda }\ |A_{i}\cap \lambda
|=\lim_{\lambda \uparrow \Lambda }\frac{\left\vert \mathbb{N}^{+}\cap
\lambda \right\vert }{p}=\frac{\mathbf{\alpha }}{p}.
\end{equation*}

(2) - We have that 
\begin{eqnarray*}
\left\vert E_{p}\cap \lambda \right\vert &=&\left\vert \{x\in \mathbb{N}%
^{+}\mid x^{p}\leq m!^{m!}\}\right\vert =\left\vert \{x\in \mathbb{N}%
^{+}\mid x\leq m!^{\frac{m!}{p}}\}\right\vert \\
&=&m!^{\frac{m!}{p}}=\left\vert \mathbb{N}^{+}\cap \lambda \right\vert ^{%
\frac{1}{p}}
\end{eqnarray*}%
Then%
\begin{equation*}
\mathfrak{num}(E_{p})=\lim_{\lambda \uparrow \Lambda }\ \left\vert \mathbb{N}%
^{+}\cap \lambda \right\vert ^{\frac{1}{p}}=\mathbf{\alpha }^{\frac{1}{p}}
\end{equation*}

$\square $

\begin{remark}
As we already observed, the choice of $\mathfrak{S}(\mathbb{N})$ is not
intrinsic, and it has been done so to make it possible to have the
properties listed in Proposition \ref{basic}. Some additional motivations
for this choice of $\mathfrak{S}(\mathbb{N})$ can be found in \cite{BDN2018}%
. Different algebraic properties of the numerosity can be rather easily
obtained by changing the choice of $\Pi $. Or if you like, you can assume
that there is an other copy of the natural numbers $\mathbf{N}\subset 
\mathbf{Ato,}$ $\mathbb{N\cap }\mathbf{N}=\varnothing ,$ and a labelling $%
\mathfrak{S}\left( \mathbf{N}\right) :=\left\{ \{0,...,n\}\ |\ n\in \Pi
_{1}\right\} .$
\end{remark}

\subsubsection{Numerosity of the rational numbers\label{raz}}

In order to get reasonable properties for the sets of rational numbers, we
set%
\begin{equation*}
\mathfrak{s}(n):=\left\{ \frac{s}{n}\mid s\in \mathbb{Z},\ -n^{2}\leq
s<n^{2}\right\}
\end{equation*}%
and%
\begin{equation*}
\mathfrak{S}(\mathbb{Q}):=\left\{ \mathfrak{s}(n)\ |\ n\in \Pi \right\} .
\end{equation*}%
Clearly, $\mathfrak{S}\left( \mathbb{Q}\right) \ $is compatible with$\ 
\mathfrak{S}\left( \mathbb{N}\right) $ and hence we can take the set $%
\mathfrak{S}\left( \mathbb{Q}\right) \cup \mathfrak{S}\left( \mathbb{N}%
\right) =\mathfrak{S}\left( \mathbb{Q}\right) $ and, via Th. \ref{laura} the
induced lattice-tree $\left( \mathfrak{L}_{\mathbb{Q}},\subseteq \right) .$

For every $q\in \mathbb{Q}$, it holds%
\begin{equation*}
\ell _{_{\mathbb{N}}}(q)=\ell _{_{\mathbb{N}}}(\{q\})=\mathfrak{s}(n)
\end{equation*}%
where $n=\min \{n\in \mathbb{N}\ |\ q\in \mathfrak{s}(n)\}.$ This choice of $%
\mathfrak{S}(\mathbb{Q}),$ allows to get the following result:

\begin{theorem}
\label{n raz} The following properties hold:

\begin{enumerate}
\item $\mathfrak{num}\left( \mathbb{Q}\cap (0,1]\right) =\mathbf{\alpha }$;

\item for all $p,q\in \mathbb{Q}$ with $p<q$, $\mathfrak{num}\left( \mathbb{Q%
}\cap (p,q]\right) =\left( p-q\right) \mathbf{\alpha }$;

\item for all $p,q\in \mathbb{R}$ with $p<q$, 
\begin{equation*}
\frac{\mathfrak{num}\left( \mathbb{Q}\cap (p,q]\right) }{\mathbf{\alpha }}%
\sim \left( p-q\right) ;
\end{equation*}

\item {}$\mathfrak{num}\left( \mathbb{Q}^{+}\right) =\mathbf{\alpha }^{2};$

\item $\mathfrak{num}\left( \mathbb{Q}\right) =2\mathbf{\alpha }^{2}+1;$

\item if $E\subset \mathbb{Q}$ is a bounded set, then, $\forall q\in \mathbb{%
Q}$ 
\begin{equation*}
\mathfrak{num}\left( q+E\right) =\mathfrak{num}\left( E\right) .
\end{equation*}
\end{enumerate}
\end{theorem}

\textbf{Proof:} (1) Take $\lambda \in C\left[ \ell _{_{\mathbb{Q}}}(n)\right]
;\ $we have that 
\begin{equation*}
\left\vert \left( \mathbb{Q}\cap (0,1]\right) \cap \mathbb{\lambda }%
\right\vert =n=\left\vert \mathbb{N}^{+}\cap \mathbb{\lambda }\right\vert ;
\end{equation*}%
the conclusion follows by taking the $\Lambda $-limit.

(2) - If $p=\frac{s_{1}}{n_{1}},\ q=\frac{s_{2}}{n_{2}}$, take $\lambda \in C%
\left[ \mathfrak{s}(m)\right] $ with $m$ larger than $|p|,|q|,n_{1}n_{2}$.
Then 
\begin{equation*}
\left( \mathbb{Q}\cap \lambda \right) =(p-q)m=(p-q)|\mathbb{N}^{+}\cap 
\mathbb{\lambda }|;
\end{equation*}%
the conclusion follows by taking the $\Lambda $-limit.

(3) - Take $\varepsilon \in \mathbb{R}^{+}$ and four numbers $p_{\varepsilon
}^{\pm },q_{\varepsilon }^{\pm }\in \mathbb{R}$ such that%
\begin{equation*}
p_{\varepsilon }^{-}<p<p_{\varepsilon }^{+}<q_{\varepsilon
}^{-}<q<q_{\varepsilon }^{+};\ \ p_{\varepsilon }^{+}-p_{\varepsilon
}^{-}\leq \varepsilon ,\ \ q_{\varepsilon }^{+}-q_{\varepsilon }^{-}\leq
\varepsilon ;
\end{equation*}%
then, 
\begin{equation*}
\mathfrak{num}\left( (p_{\varepsilon }^{+},q_{\varepsilon }^{-}]\right) \leq 
\mathfrak{num}\left( (p,q]\right) \leq \mathfrak{num}\left( (p_{\varepsilon
}^{-},q_{\varepsilon }^{+}]\right)
\end{equation*}%
and by (2) 
\begin{equation*}
\left[ \left( p-q\right) -2\varepsilon \right] \mathbf{\alpha }\leq \ 
\mathfrak{num}\left( (p,q]\right) \leq \left[ \left( p-q\right)
-2\varepsilon \right] \mathbf{\alpha ;}
\end{equation*}%
hence,%
\begin{equation*}
\left\vert \frac{\ \mathfrak{num}\left( (p,q]\right) }{\mathbf{\alpha }}%
-\left( p-q\right) \right\vert \leq 2\varepsilon
\end{equation*}

(4) - The map $\Phi :(0,1]\times \mathbb{N}^{+}\rightarrow \mathbb{Q}^{+}$
defined by 
\begin{equation*}
\Phi (x,n)=(n-1]+x
\end{equation*}%
is a comparison map, since, for a suitable $m\in \mathbb{N},$ 
\begin{equation*}
\ell (\Phi (x,n))=\ell \left( (n-1]+x\right) =\mathfrak{s}(m)=\ell \left(
(x,n)\right)
\end{equation*}%
and hence, by Cor. \ref{20}, 
\begin{eqnarray*}
\mathfrak{num}\left( \mathbb{Q}^{+}\right) &=&\mathfrak{num}\left(
(0,1]\times \mathbb{N}^{+}\right) =\mathfrak{num}\left( [0,1)\right) \cdot 
\mathfrak{num}\left( \mathbb{N}^{+}\right) \\
&=&\mathfrak{num}\left( [0,1)\right) \cdot \mathfrak{num}\left( \mathbb{N}%
^{+}\right) =\mathbf{\alpha }^{2}.
\end{eqnarray*}

(5) - Since $\left\vert \mathbb{Q}^{+}\cap \lambda \right\vert =\left\vert 
\mathbb{Q}^{-}\cap \lambda \right\vert ,$ we have that $\mathfrak{num}\left( 
\mathbb{Q}^{+}\right) =\mathfrak{num}\left( \mathbb{Q}^{-}\right) .$ Then%
\begin{equation*}
\mathfrak{num}\left( \mathbb{Q}^{+}\right) =\mathfrak{num}\left( \mathbb{Q}%
^{+}\right) +\mathfrak{num}\left( \mathbb{Q}^{-}\right) +\mathfrak{num}%
\left( \{0\}\right) =2\mathbf{\alpha }^{2}+1.
\end{equation*}

(6) - If $E\subset \mathfrak{s}(n_{1})$ and $q\leq n_{2},$ then there exists 
$n$ sufficiently large that $E\subset \mathfrak{s}(n)$ and $q+E\subset 
\mathfrak{s}(n);$ then $\forall \lambda \in C\left[ \mathfrak{s}(n)\right] ,$%
we have that 
\begin{equation*}
\left\vert E\cap \lambda \right\vert =\left\vert \left( q+E\right) \cap
\lambda \right\vert
\end{equation*}

$\square $

\subsubsection{Numerosity of the real numbers\label{NR}}

For every $\ \Xi \in \wp _{\omega }(\left( 0,1\right] )$ and every $n\in \Pi
,$ we set%
\begin{equation*}
\mathfrak{s}\left( \Xi ,n\right) =\dbigcup\limits_{q\in \mathfrak{s}%
(n)}\left( q+\frac{1}{n}\Xi \right) ,\ n\in \Pi
\end{equation*}%
\begin{equation*}
\mathfrak{S}(\mathbb{R}):=\left\{ \mathfrak{s}\left( \Xi ,n\right) \ |\ \Xi
\in \wp _{\omega }(\left( 0,1\right] ),\ n\in \Pi \right\} ;
\end{equation*}%
hence an element $\mathfrak{s}\in \mathfrak{S}(\mathbb{R})$ is contained in $%
[-n,n)$ and is the union of $2n^{2}$ copies of $\frac{1}{n}\Xi =\left\{ 
\frac{r}{n}\ |\ r\in \Xi \right\} $, each of them contained in an interval $%
[q,\ q+\frac{1}{n})$. Clearly, $\mathfrak{S}\left( \mathbb{R}\right) \ $and $%
\mathfrak{S}\left( \mathbb{Q}\right) $ are compatible\ and once again, we
apply Th. \ref{laura} to the lattice-three $(\mathfrak{S}\left( \mathbb{R}%
\right) \cup \mathfrak{S}\left( \mathbb{Q}\right) ,\subseteq )=(\mathfrak{S}%
\left( \mathbb{R}\right) ,\subseteq )\ $and we get the induced label-tree $%
\left( \mathfrak{L}_{\mathbb{R}},\subseteq \right) $.

Now, we need to give a name to $\mathfrak{num}\left( (0,1]\right) ;$
inspired by (\ref{alfa0}) and Prop. \ref{n raz}-(i), we put 
\begin{equation*}
\mathbf{\beta }:=\mathfrak{num}\left( (0,1]\right)
\end{equation*}

\begin{theorem}
\label{n reali} The following properties hold

\begin{enumerate}
\item for all $n\in \mathbb{N}$, $\mathfrak{num}\left( (n,n+1]\right) =%
\mathbf{\beta }$;

\item for all $p,q\in \mathbb{Q}$ with $p<q$, $\mathfrak{num}\left( \left[
p,q\right) \right) =\left( p-q\right) \mathbf{\beta }$;

\item for all $p,q\in \mathbb{R}$ with $p<q$, $\mathfrak{num}\left( \left[
p,q\right) \right) \sim \left( p-q\right) \mathbf{\beta }$;

\item $\mathfrak{num}\left( \mathbb{R}^{+}\right) =\mathbf{\alpha }\mathbf{%
\beta }$;

\item $\mathfrak{num}\left( \mathbb{R}\right) =2\mathbf{\alpha }\mathbf{%
\beta }+1$;

\item if $E\subset \mathbb{E}$ is a bounded set, then, $\forall q\in \mathbb{%
Q}$ 
\begin{equation*}
\mathfrak{num}\left( q+E\right) =\mathfrak{num}\left( E\right) .
\end{equation*}
\end{enumerate}
\end{theorem}

\textbf{Proof: }It is similar to the proof of Prop. \ref{n raz} with minor
changes.

$\square $

\bigskip

It is interesting to compare the numerosity of measurable subsets of $%
\mathbb{R}$ with their Lebesgue measure.

\begin{definition}
For every $A\in \Lambda $, and every $\gamma \in \mathcal{N}$, we set 
\begin{equation*}
m_{\gamma }(A):=st\left( \frac{\mathfrak{num}(A)}{\gamma }\right)
\end{equation*}%
We will call $m_{\gamma }(A)$ $\gamma $-measure of $A.$
\end{definition}

\begin{lemma}
\label{sc}The $\gamma $-measure satisfies the following properties:

\begin{enumerate}
\item it is finitely additive: for all sets $A,B$ 
\begin{equation*}
m_{\gamma }\left( A\cup B\right) =m_{\gamma }\left( A\right) +m_{\gamma
}\left( B\right) -m_{\gamma }\left( A\cap B\right) ;
\end{equation*}

\item it is superadditive, namely given a denumerable partition $\left\{
A_{n}\right\} _{n\in \mathbb{N}}$ of a set $A\subset \mathbb{R}$, then 
\begin{equation*}
m_{\gamma }\left( A\right) \geq \sum_{n=0}^{\infty }m_{\gamma }\left(
A_{n}\right) .
\end{equation*}
\end{enumerate}
\end{lemma}

\textbf{Proof -} (i) This is a trivial consequence of the additivity of the
numerosity.

(ii) We have that for all $N\in \mathbb{N}$, 
\begin{equation*}
\mathfrak{num}\left( A\right) \geq \mathfrak{num}\left(
\dbigcup\limits_{n=0}^{N}A_{n}\right) =\sum_{n=0}^{N}\mathfrak{num}\left(
A_{n}\right) ,
\end{equation*}%
hence 
\begin{equation*}
st\left( \frac{\mathfrak{num}\left( A\right) }{\gamma }\right) \geq st\left(
\sum_{n=0}^{N}\frac{\mathfrak{num}\left( A_{n}\right) }{\gamma }\right)
=\sum_{n=0}^{N}st\left( \frac{\mathfrak{num}\left( A_{n}\right) }{\gamma }%
\right) ;
\end{equation*}%
therefore, 
\begin{equation*}
m_{\gamma }\left( A\right) \geq \sum_{n=0}^{N}m_{\gamma }\left( A_{n}\right)
.
\end{equation*}%
The conclusion follows taking the Cauchy limit in the above inequality for $%
N\rightarrow \infty $.

$\square $

\begin{definition}
A set $A\subset \mathbb{R}^{N}\ $such that $m_{\gamma }(A)<\infty ,$ is
called $m_{\gamma }$-measurable if 
\begin{equation*}
m_{\gamma }(A)=\sup \left\{ m_{\gamma }\left( F\right) \ |\ F\subset A,\ F\ 
\text{is\ closed}\right\}
\end{equation*}
\end{definition}

\begin{lemma}
The family $\mathfrak{M}_{\gamma }^{<\infty }$ of $m_{\gamma }$-measurable
sets form a $\sigma $-ring.
\end{lemma}

\textbf{Proof}: First let us prove that $\mathfrak{M}_{\gamma }^{<\infty }$
is a ring. Let $A,B\in \mathfrak{M}_{\gamma }^{<\infty };$ we have show that 
$A\cap B\in \mathfrak{M}_{\gamma }^{<\infty }$. Take $\varepsilon \in 
\mathbb{R}^{+}$ and two sets $F\subset A$ and $G\subset B$ such that $%
m_{\gamma }\left( A\right) -m_{\gamma }\left( F\right) \leq \varepsilon /2,\
m_{\gamma }\left( B\right) -m_{\gamma }\left( G\right) \leq \varepsilon /2.$
Then%
\begin{eqnarray*}
m_{\gamma }\left( A\cap B\right) -m_{\gamma }\left( F\cap G\right)
&=&m_{\gamma }\left( \left( A\cap B\right) \backslash \left( F\cap G\right)
\right) \\
&\leq &m_{\gamma }\left( \left[ A\backslash F\right] \cup \left[ B\backslash
G\right] \right) \\
&\leq &m_{\gamma }\left( A\backslash F\right) +m_{\gamma }\left( B\backslash
G\right) \leq \varepsilon
\end{eqnarray*}%
By the arbitrariness of $\varepsilon $ the conclusion follows. Now let us
prove the $\sigma $-additivity. Let $\left\{ A_{n}\right\} _{n\in \mathbb{N}%
} $ be a denumerable partition $\left\{ A_{n}\right\} _{n\in \mathbb{N}}$ of
a set $A$ and let $\varepsilon \in \mathbb{R}^{+};$ since $%
\sum_{n=0}^{\infty }m_{\gamma }\left( A_{n}\right) \leq m_{\gamma }\left(
A\right) $, this series is convergent and hence there exists $m$ such that $%
\sum_{n=m+1}^{\infty }m_{\gamma }\left( A_{n}\right) \leq \varepsilon /2;$
now take a family of closed sets $F_{n}$ such that $m_{\gamma }\left(
A_{n}\right) -m_{\gamma }\left( F_{n}\right) \leq 2^{-(n+2)}\varepsilon .$
Then, putting $F=\dbigcup\limits_{n=0}^{m}F_{n}$ we have that $F$ is closed
and 
\begin{eqnarray*}
m_{\gamma }\left( A\right) -m_{\gamma }\left( F\right) &=&m_{\gamma }\left(
\dbigcup\limits_{n=0}^{\infty }A_{n}\right) -m_{\gamma }\left(
\dbigcup\limits_{n=0}^{m}F_{n}\right) \\
&\leq &\sum_{n=0}^{\infty }m_{\gamma }\left( A_{n}\right)
-\sum_{n=0}^{m}m_{\gamma }\left( F_{n}\right) \\
&=&\sum_{n=0}^{m}\left[ m_{\gamma }\left( A_{n}\right) -m_{\gamma }\left(
F_{n}\right) \right] +\sum_{n=m+1}^{\infty }m_{\gamma }\left( A_{n}\right) \\
&\leq &\sum_{n=0}^{m}2^{-(n+2)}\varepsilon -\frac{\varepsilon }{2}%
<\varepsilon .
\end{eqnarray*}

$\square $

From this lemma, the following theorem immediately follows:

\begin{theorem}
If $\mathfrak{M}_{\gamma }$ is the $\sigma $-algebra generated by $\mathfrak{%
M}_{\gamma }^{<\infty };\ $then $m_{\gamma }:\mathfrak{M}_{\gamma
}\rightarrow \left[ 0,\infty \right] $ is a measure.
\end{theorem}

\bigskip

By the above theorem and well known results of measure theory we get the
following facts:

\begin{itemize}
\item the Borellian sets are $m_{\gamma }$-measurable;

\item if $A$ is $m_{\gamma }$-measurable, then%
\begin{equation*}
m_{\gamma }(A)=\inf \left\{ m_{\gamma }\left( H\right) \ |\ A\subset H,\ F\ 
\text{is\ open}\right\} =\sup \left\{ m_{\gamma }\left( F\right) \ |\
F\subset A,\ F\ \text{is\ closed}\right\}
\end{equation*}
\end{itemize}

From here, it follows that

\begin{corollary}
If we take $\gamma =\beta ,$ then $m_{\beta }$ is equal to the Lebesgue
measure $m_{L}$.
\end{corollary}

\begin{corollary}
If $A\subset \mathbb{R}$ is a Lebesgue measurable set, then 
\begin{equation*}
num(A)=(1+\varepsilon )\mathbf{\beta \,}m_{L}\left( A\right) ,\ \
\varepsilon \sim 0.
\end{equation*}
\end{corollary}

\begin{corollary}
If $A\subset \mathbb{R}^{N}$ is a Lebesgue measurable set, then 
\begin{equation*}
m_{N}\left( A\right) =st\left( \frac{\mathfrak{num}(A)}{\mathbf{\beta }^{N}}%
\right) .
\end{equation*}%
where $m_{N}$ denotes the $N$-dimensional Lebesgue measure.
\end{corollary}

\textbf{Proof}: By Th. \ref{n reali}, we have that for all $p,q\in \mathbb{Q}
$ with $p<q$, $\mathfrak{num}\left( \left[ p,q\right) \right) =\left(
p-q\right) \mathbf{\beta }$; hence 
\begin{equation*}
\mathfrak{num}\left( \left[ p_{1},q_{1}\right) \times ...\times \left[
p_{N},q_{N}\right) \right) =\left( q_{1}-p_{1}\right) \cdot ....\cdot \left(
q_{N}-p_{N}\right) \mathbf{\beta }^{N}
\end{equation*}%
Then, if we take the $\gamma $-measure with $\gamma =\mathbf{\beta }^{N},$
on the $\sigma $-algebra of the $\mathbf{\beta }^{N}$-measurable sets, we
have that 
\begin{equation*}
m_{N}\left( A\right) =m_{\mathbf{\beta }^{N}}\left( A\right) =st\left( \frac{%
\mathfrak{num}(A)}{\mathbf{\beta }^{N}}\right) .
\end{equation*}

$\square $

\bigskip

\subsection{Complementary examples}

In this section, we will give some examples to show how some special
properties of the numerosities can be implemented in different circumstances.

\subsubsection{Exponentiation of ordinal numerosities\label{EON}}

In this section, we will examine the relation of the ordinal exponentiation $%
\beta ^{\left\langle \gamma \right\rangle }$ and the numerosity
exponentiation $\beta ^{\gamma }$. It is not possible to require that $%
\forall \gamma \in \mathbf{Ord},\ \beta ^{\left\langle \gamma \right\rangle
}=\beta ^{\gamma }$ since the map $\left\langle \gamma \right\rangle \mapsto
\beta ^{\left\langle \gamma \right\rangle }$ has fixed points. For example,
for every $n\in \mathbb{N},$ we have that $n^{\left\langle \omega
\right\rangle }=\omega \neq n^{\omega }.$ However, it is natural to
investigate when an ordinal number written in the Cantor normal form $%
\dsum\limits_{k=0}^{n}b_{k}\omega ^{\left\langle j_{k}\right\rangle }$ is
equal to the numerosity $\dsum\limits_{k=0}^{n}b_{k}\omega ^{j_{k}}.$ The
fix points of the ordinal exponentiation $\omega ^{\left\langle \omega
\right\rangle }$ are called $\varepsilon $-numbers and they are denoted
denote by the symbol $\varepsilon _{j},\ j\in \mathbf{Ord}$. In particular
the smallest of them is given by:%
\begin{equation*}
\varepsilon _{0}:=\sup_{\mathbf{Ord}}\left\{ \ \omega ^{\left\langle \omega
^{\left\langle \omega ^{...}\right\rangle }\right\rangle }\right\} .
\end{equation*}

Then, $\omega ^{\varepsilon _{0}}>\varepsilon _{0}=\omega ^{\left\langle
\varepsilon _{0}\right\rangle }$ and, by theorem \ref{lana}, it follows
that, if $\beta \geq \varepsilon _{0},$ 
\begin{equation*}
\beta =\dsum\limits_{k=0}^{n}b_{k}\omega ^{\left\langle j_{k}\right\rangle
}\,<\dsum\limits_{k=0}^{n}b_{k}\omega ^{j_{k}}.
\end{equation*}%
However, it is natural to require that%
\begin{equation*}
\dsum\limits_{k=0}^{n}b_{k}\omega ^{\left\langle j_{k}\right\rangle
}\,=\dsum\limits_{k=0}^{n}b_{k}\omega ^{j_{k}}\ \ \text{when}\ \ \beta
<\varepsilon _{0}
\end{equation*}%
This request can be satisfied if we choose a suitable set $\mathfrak{S}_{%
\mathbf{Ord}}$ and the induced label tree. To this aim, it is convenient to
set 
\begin{equation*}
\mathfrak{s}(\beta ):=\left\{ \dsum\limits_{k=0}^{n}c_{k}\omega
^{\left\langle j_{k}\right\rangle }\mid c_{k}\sqsubseteq _{\mathbb{N}%
}b_{k}\right\} ,\ \beta =\dsum\limits_{k=0}^{n}b_{k}\omega
\end{equation*}%
and%
\begin{equation*}
\mathfrak{S}(\mathbf{Ord}):=\left\{ \mathfrak{s}(\beta )\ |\ \beta \in 
\mathbf{Ord}\right\} .
\end{equation*}%
The set $\mathfrak{S}(\mathbf{Ord})\ $is compatibel with $\mathfrak{S}(%
\mathbb{R)}$ since $\mathfrak{S}(\mathbf{Ord})\cap \mathfrak{S}(\mathbb{R)}=%
\mathfrak{S}(\mathbb{N)}\subset \mathfrak{S}(\mathbb{R)}$. Then we put $%
\mathfrak{S}(\mathbf{O}):=\mathfrak{S}(\mathbf{Ord})\cup \ \mathfrak{S}(%
\mathbb{R)}$ and we denote with $\left( \mathfrak{L}_{\mathbf{O}},\subseteq
\right) $ the label-tree induced by $\mathfrak{S}(\mathbf{O}).$

Now let us check some properties of the label $\ell _{\mathbf{O}}(a).$

\begin{lemma}
\label{ven}The label $\ell _{\mathbf{O}}(\gamma )$ of a ordinal number
satisfies the following properties:

\begin{itemize}
\item (i) $\ell _{\mathbf{O}}(\omega ^{\left\langle j\right\rangle
})=\{\omega ^{\left\langle j\right\rangle }\}$

\item (ii) $\ell _{\mathbf{O}}\left( c\omega ^{\left\langle j\right\rangle
}\right) =\ell _{\mathbb{N}}\left( c\right) \vee \ell _{\mathbf{O}}(\omega
^{\left\langle j\right\rangle });$

\item (iii) if $\gamma =\dsum\limits_{k=0}^{n}b_{k}\omega ^{\left\langle
j_{k}\right\rangle },\ $then $\ \ell _{\mathbf{O}}(\gamma
)=\dbigvee\limits_{k=0}^{n}\left[ \ell _{\mathbb{N}}\left( c_{k}\right) \vee
\ell _{\mathbf{O}}(\omega ^{\left\langle j\right\rangle })\right] \ $
\end{itemize}
\end{lemma}

\textbf{Proof}: - (i) follows from the fact that the $\omega ^{\left\langle
j\right\rangle }$'s have no predecessor. Also we have that 
\begin{eqnarray*}
\ell _{\mathbf{O}}\left( c\omega ^{\left\langle j\right\rangle }\right)
&=&\left\{ x\in \mathbf{Ord}\ |\ x\sqsubseteq _{\mathbf{O}}c\omega
^{\left\langle j\right\rangle }\right\} \\
&=&\left\{ b\omega ^{\left\langle j\right\rangle }\ |\ b\sqsubseteq _{%
\mathbb{N}}c\right\} =\ell _{\mathbb{N}}\left( c\right) \vee \ell _{\mathbf{O%
}}(\omega ^{\left\langle j\right\rangle }).
\end{eqnarray*}

(iii) - If $\gamma \in \mathbf{O}\left( \Theta _{1}\right)
=\dsum\limits_{k=0}^{n}b_{k}\omega ^{\left\langle j_{k}\right\rangle },$ 
\begin{eqnarray*}
\ell _{\mathbf{O}}\left( \gamma \right) &=&\mathfrak{s}(\gamma ):=\left\{
\dsum\limits_{k=0}^{n}b_{k}\omega ^{\left\langle j_{k}\right\rangle }\mid
b_{k}\omega ^{\left\langle j_{k}\right\rangle }\sqsubseteq _{\mathbb{N}%
}c_{k}\omega ^{\left\langle j_{k}\right\rangle }\right\} \\
&=&\dbigvee\limits_{k=0}^{n}\left[ \ell _{\mathbb{N}}\left( c_{k}\right)
\vee \ell _{\mathbf{O}}(\omega ^{\left\langle j\right\rangle })\right] .
\end{eqnarray*}

$\square $

\bigskip

The next theorem characterizes the numerosity exponentiation between
ordinals.

\begin{theorem}
\label{TTT}If $\gamma \in \mathbf{Ord,}$ then%
\begin{equation*}
\omega ^{\gamma }=\mathfrak{num}\left( \mathfrak{F}_{fin}\left( \mathbf{O(}%
\gamma \mathbf{)},\mathbb{N}^{+}\right) \right)
\end{equation*}
\end{theorem}

\textbf{Proof}: Let us consider the map $\Psi :\mathfrak{F}_{fin}\left( 
\mathbf{O(}\gamma \mathbf{)},\mathbb{N}^{+}\right) \rightarrow \mathbf{Ord}$
defined as follows%
\begin{equation*}
\Psi \left( f\right) =\sum_{\xi \in D_{f}}f\left( \xi \right) \omega
^{\left\langle \xi \right\rangle }.
\end{equation*}%
where $D_{f}\in \wp _{\omega }(\mathbf{O(}\gamma \mathbf{)})$ is the domain
of $f$. By Lemma \ref{ven}-(\textit{iii}), we get%
\begin{equation*}
\ell _{\mathbf{O}}\left( \Psi (f)\right) =\dbigvee\limits_{\xi \in D_{f}}%
\left[ \ell _{\mathbf{O}}(\xi )\vee \ell _{\mathbf{O}}(f\left( \xi \right) )%
\right] .
\end{equation*}%
Since the function $f$ is identified with its graph which is a finite set,
by Prop.\ref{PP}-(q4-q1), we have that%
\begin{eqnarray*}
\ell _{\mathbf{O}}(f) &=&\ell _{\mathbf{O}}\left( \left\{ (\xi ,f\left( \xi
\right) )\ |\ \xi \in D_{f}\right\} \right) \\
&=&\dbigvee\limits_{\xi \in D_{f}}\ell _{\mathbf{O}}(\left( \xi ,f\left( \xi
\right) \right) )=\dbigvee\limits_{\xi \in D_{f}}\left[ \ell _{\mathbf{O}%
}(\{\{\xi \}\})\vee \ell _{\mathbf{O}}(\{\{f\left( \xi \right) \}\})\right]
\end{eqnarray*}%
Then 
\begin{eqnarray*}
\ell _{\mathbf{O}}(f)\cap \mathbf{Ord} &=&\left( \dbigvee\limits_{\xi \in
D_{f}}\left[ \ell _{\mathbf{O}}(\{\{\xi \}\})\vee \ell _{\mathbf{O}%
}(\{\{f\left( \xi \right) \}\})\right] \right) \cap \mathbf{Ord} \\
&=&\dbigvee\limits_{\xi \in D_{f}}\left( \left[ \ell _{\mathbf{O}}(\{\{\xi
\}\})\cap \mathbf{Ord}\right] \vee \left[ \ell _{\mathbf{O}}(\{\{f\left( \xi
\right) \}\})\cap \mathbf{Ord}\right] \right) \\
&=&\dbigvee\limits_{\xi \in D_{f}}\left[ \ell _{\mathbf{O}}(\xi )\vee \ell _{%
\mathbf{O}}(f\left( \xi \right) )\right] =\ell _{\mathbf{O}}\left( \Psi
(f)\right)
\end{eqnarray*}%
Hence, by Prop.\ref{strega}, $\Psi $ is a comparison bijection and by Prop.%
\ref{annamaria}%
\begin{equation*}
\mathfrak{num}\left( \mathfrak{F}_{fin}\left( \mathbf{O(}\gamma \mathbf{)},%
\mathbb{N}^{+}\right) \right) =\mathfrak{num}\left( \mathbb{N}\right) ^{%
\mathfrak{num}\left( \mathbf{O(}\gamma \mathbf{)}\right) }=\omega ^{\gamma }
\end{equation*}

$\square $

\begin{theorem}
\label{TTTT}For every $\gamma \in \mathbf{Ord}$, if $\gamma <\varepsilon
_{0} $, 
\begin{equation}
\omega ^{\gamma }=\omega ^{\left\langle \gamma \right\rangle }  \label{aveM}
\end{equation}
\end{theorem}

\textbf{Proof}: We argue by induction over $\gamma $. If $\gamma =0,$ (\ref%
{aveM}) holds trivially, and it holds also if $\gamma $ is a successor. Now
let us assume that (\ref{aveM}) holds $\forall \beta <\gamma $ and let us
prove it for $\gamma :$%
\begin{equation*}
\omega ^{\gamma }=\mathfrak{num}\left( \mathfrak{F}_{fin}\left( \mathbf{O(}%
\gamma \mathbf{)},\mathbb{N}^{+}\right) \right) =\mathfrak{num}\left(
\dbigcup\limits_{\beta <\gamma }\mathfrak{F}_{fin}\left( \mathbf{O(}\beta 
\mathbf{)},\mathbb{N}^{+}\right) \right) \leq \min \{\tau \in \mathbf{Ord}\
|\ \tau >\omega ^{\beta }\}
\end{equation*}%
Since $\tau <\varepsilon _{0},$ 
\begin{equation*}
\omega ^{\gamma }\leq \min \{\tau \in \mathbf{Ord}\ |\ \tau >\omega ^{\beta
}\}\leq \underset{\mathbf{Ord}}{\sup }\{\omega ^{\beta }\ |\ \beta <\gamma \}
\end{equation*}%
By the inductive assumption, if $\beta <\gamma ,$ then $\omega ^{\beta
}=\omega ^{\left\langle \beta \right\rangle }$ and hence%
\begin{equation*}
\omega ^{\gamma }\leq \underset{\mathbf{Ord}}{\sup }\{\omega ^{\left\langle
\beta \right\rangle }\ |\ \beta <\gamma \}=\omega ^{\left\langle \gamma
\right\rangle }.
\end{equation*}%
The conclusion follows from Th.\ref{CAC}.

$\square $

In conclusion, the set $\mathbf{O}(\varepsilon _{0})$ is closed for
exponentiation and the natural ordinal operations and the numerosity
operation coincide.

\subsubsection{$\beth _{1}$ versus $\mathbf{\protect\beta }$\label{bb}}

Probably the first set having the cardinality of continuum which comes to
your mind is either $[0,1]$ or $\wp (\mathbb{N});$ we have seen that $%
\mathfrak{num}([0,1])=\mathbf{\beta }+1$ and $\mathfrak{num}(\wp (\mathbb{N}%
))=\beth _{1}.$ It is natural to establish a relation between them. This can
be done by choosing a suitable comparison map. Probably the most natural way
to map $\wp (\mathbb{N})$ over $\left[ 0,1\right] $ is the binary expansion
of a real number given by%
\begin{equation*}
\psi (B):=\dsum\limits_{n\in B}2^{-\left( n+1\right) }.
\end{equation*}%
Hence, if $B\subset \wp (\mathbb{N})$ contains the number $n,$ then, the $%
(n+1)$-th digit of the dual expansion of $\psi (B)$ is "1". $\psi $ cannot
be a comparison map since it is not injective. However, its restriction to
infinite set 
\begin{equation*}
\psi _{\text{\textsc{r}}}:\mathbf{Inf}(\mathbb{N})\rightarrow (0,1],\ \ 
\mathbf{Inf}(\mathbb{N}):=\wp (\mathbb{N})\cap \mathbf{Inf}
\end{equation*}%
is bijective. If we want $\psi _{\text{\textsc{r}}}$ to be a comparison map,
it is sufficient to introduce the set 
\begin{equation*}
\mathfrak{S}(\psi _{\text{\textsc{r}}})=\{\{x,\psi _{\text{\textsc{r}}%
}(x)\}\ |\ x\in \mathbf{Inf}(\mathbb{N})\};
\end{equation*}%
The sets $\mathfrak{S}(\psi _{\text{\textsc{r}}})$ and $\mathfrak{S}(\mathbf{%
O)}$ are trivially compatible since $\mathfrak{S}(\psi _{\text{\textsc{r}}%
})\cap \mathfrak{S}(\mathbf{O})=\varnothing $ We will denote by $(\mathfrak{L%
}_{\mathbf{O,}\psi _{\text{\textsc{r}}}},\subseteq )$ the inducede label
tree. Using the labelling $\mathfrak{L}_{\mathbf{O,}\psi _{\text{\textsc{r}}%
}},\ $we get the following results:

\begin{theorem}
The numerosity of the unit interval is given by%
\begin{equation*}
\mathfrak{num}(\left[ 0,1\right] )=\beth _{1}-2^{\omega }+1
\end{equation*}
\end{theorem}

\textbf{Proof}: By our construction $\psi _{\text{\textsc{r}}}$ is a
comparison map, then applying the rules of numerosity:%
\begin{eqnarray*}
\mathfrak{num}(\left[ 0,1\right] ) &=&\mathfrak{num}((0,1])+1=\mathfrak{num}%
(\psi _{\text{\textsc{r}}}^{-1}\left( (0,1]\right) )+1 \\
&=&\mathfrak{num}(\mathbf{Inf}\left( \mathbb{N}\right) )+1=\mathfrak{num}%
(\wp (\mathbb{N})\backslash \wp _{\omega }(\mathbb{N}))+1 \\
&=&\mathfrak{num}(\wp (\mathbb{N}))-\mathfrak{num}(\wp _{\omega }(\mathbb{N}%
))+1 \\
&=&\beth _{1}-2^{\mathfrak{\omega }}+1.
\end{eqnarray*}

$\square $

\begin{corollary}
We have that%
\begin{equation*}
\mathbf{\beta }=\beth _{1}-2^{\omega }
\end{equation*}
\end{corollary}

\end{document}